\documentclass[sn-basic]{sn-jnl}
\usepackage[utf8]{inputenc}
\usepackage[T1]{fontenc}
\usepackage[title]{appendix}
\usepackage{amsmath,amsthm,amssymb,amsfonts,mathtools}
\usepackage{xspace}
\usepackage{makecell}
\usepackage{enumitem}
\usepackage[authoryear]{natbib}
\usepackage{tabularx}
\usepackage{array,multirow,graphicx,adjustbox}
\usepackage{xcolor}
\usepackage{hyperref}
\usepackage{comment}
\usepackage{pgfplots}
\usepackage{enumitem}
\usepackage{booktabs}
\usepackage{siunitx}
\usepackage{algorithm}
\usepackage{longtable}
\usepackage{array}
\usepackage{booktabs}
\usepackage{xltabular}
\usepackage{array}
\usepackage{algpseudocode}
\usepackage{caption}
\pgfplotsset{compat=1.18}
\usepgfplotslibrary{groupplots}
\usepgfplotslibrary{colormaps}
\usetikzlibrary{matrix,decorations.pathreplacing,calc}

\newcommand{\NA}{\multicolumn{1}{c}{--}}

\makeatletter
\newcommand{\customlabel}[2]{\protected@write \@auxout {}{\string \newlabel {#1}{{#2}{\thepage}{#2}{#1}{}} }\hypertarget{#1}{#2}
}
\makeatother

\theoremstyle{definition}
\newtheorem{definition}{Definition}[section]

\newcommand{\setCurricula}{\mathcal{K}}
\newcommand{\setMandCourses}{\mathcal{C}^{\text{mand}}}
\newcommand{\setElectiveCourses}{\mathcal{C}^{\text{el}}}
\newcommand{\setPeriods}{\mathcal{P}}
\newcommand{\setPeriodsAtDay}{\mathcal{P}}
\newcommand{\setDays}{\mathcal{D}}
\newcommand{\setLecturers}{\mathcal{L}}
\newcommand{\setLecturersOfCourse}{\mathcal{L}}
\newcommand{\setSubsequentPeriods}{\mathcal{P}^{\text{+}}}
\newcommand{\setCourses}{\mathcal{C}}
\newcommand{\setOfCoursesInCurr}{\mathcal{C}}
\newcommand{\setCoursesOfLecturer}{\mathcal{C}}
\newcommand{\setLectures}{\mathcal{C}^{\text{lec}}}
\newcommand{\setExercises}{\mathcal{C}^{\text{ex}}}
\newcommand{\setTutorials}{\mathcal{C}^{\text{tut}}}
\newcommand{\setModules}{\mathcal{M}}

\newcommand{\setRooms}{\mathcal{R}}
\newcommand{\setPrefPeriodsOfCourse}{\mathcal{P}^{\text{pref}}}
\newcommand{\setRequiredRooms}{\mathcal{R}^{\text{req}}}
\newcommand{\setRoomConflicts}{\mathcal{X}^{\text{rooms}}}
\newcommand{\setRoomDaySetting}{\mathcal{R}^{\text{day}}}
\newcommand{\setLecturersCompactWeek}{\mathcal{L}^{\text{cons\_days}}}
\newcommand{\setLecturersCompactPeriods}{\mathcal{L}^{\text{cons\_p}}}
\newcommand{\setLecturersDistWeek}{\mathcal{L}^{\text{non\_cons\_days}}}
\newcommand{\setLecturersDistPeriods}{\mathcal{L}^{\text{non\_cons\_p}}}
\newcommand{\setLecturersPrefDaysPerWeek}{\mathcal{L}^{\text{ntd}}}

\newcommand{\setExcludedPeriods}{\mathcal{P}^{\text{ex}}}
\newcommand{\setCompatibleCourses}{\mathcal{C}}

\newcommand{\setCourseConflict}{\mathcal{X}}
\newcommand{\setfixPeriods}{\mathcal{P}^{\text{fix}}}
\newcommand{\setPeriodGroups}{\mathcal{P}^{\text{group}}}

\newcommand{\setExcludedPeriodsOfRoom}{\mathcal{P}^{\text{ex}}}
\newcommand{\setPeriodsOfLastYear}{\mathcal{P}^{\text{prev}}}
\newcommand{\setCourseTuples}{\mathcal{X}^{\text{soft}}}
\newcommand{\setSolutions}{X}
\newcommand{\setImage}{z(X)}
\newcommand{\setNondomImages}{Z_\text{nd}}
\newcommand{\setApproxSet}{Z}
\newcommand{\setApproxSolutions}{S}

\newcommand{\parPeriodPenalty}{\mathsf{pen\_stud}}

\newcommand{\parRoomPenalty}{\mathsf{pen\_room}}
\newcommand{\parNumberOfDays}{\mathsf{number\_days}}
\newcommand{\parTimesPerWeek}{\mathsf{ses}}
\newcommand{\parFollowingCourse}{\mathsf{follow}}
\newcommand{\parMaxPeriodSet}{\mathsf{max}}
\newcommand{\parCourseCurrWeight}{w}

\newcommand{\parSeatsForCourse}{\mathsf{seats\_course}}
\newcommand{\parSeatsInRoom}{\mathsf{seats\_room}}
\newcommand{\parLecturersObj}{z^{\text{lec}}}

\newcommand{\parStudentsObj}{z^{\text{stud}}}

\newcommand{\parPerturbationsObj}{z^{\text{pert}}}

\newcommand{\parStudentsBound}{\mathsf{stud\_bound}}
\newcommand{\parLecturersBound}{\mathsf{lec\_bound}}
\newcommand{\parPerturbationsBound}{\mathsf{pert\_bound}}
\newcommand{\parLecAnchorPt}{z^{\text{anc},\text{lec}}}
\newcommand{\parStudAnchorPt}{z^{\text{anc},\text{stud}}}
\newcommand{\parLecAnchorPtTheta}{z_{\text{anc}}^{\text{lec},\theta}}
\newcommand{\parStudAnchorPtTheta}{z_{\text{anc}}^{\text{stud},\theta}}

\newcommand{\parApproxNadirPt}{z^{\text{nad}}(Z)}
\newcommand{\parRefPt}{z^{\mathrm{ref}}}

\newcommand{\varAssignment}{a}
\newcommand{\varRoomAvailability}{y^{\text{room}}}
\newcommand{\varLecAtDay}{a^{\text{day}}}
\newcommand{\varLecAtDayPlusOne}{\textnormal{vio}^{\text{+1}}}
\newcommand{\varLecAtDayPlusTwo}{\textnormal{vio}^{\text{+2}}}
\newcommand{\varLecAtDayPlusThree}{\textnormal{vio}^{\text{+3}}}
\newcommand{\varLecAtDayPlusFour}{\textnormal{vio}^{\text{+4}}}
\newcommand{\varLecDeviationCoursesPerDay}{\textnormal{vio}^{\text{c\_per\_day}}}
\newcommand{\varLecAtPeriod}{a^{\text{period}}}
\newcommand{\varLecConsecutivePeriods}{\textnormal{vio}^{\text{cons}}}
\newcommand{\varLecDisjointPeriods}{\textnormal{vio}^{\text{non\_cons}}}
\newcommand{\varLecDevDaysPos}{\textnormal{vio}^{\text{pos}}}
\newcommand{\varLecDevDaysNeg}{\textnormal{vio}^{\text{neg}}}
\newcommand{\varCourseDistVio}{\textnormal{vio}^{\text{dist}}}
\newcommand{\varCurCourseVio}{\textnormal{vio}^{\text{elective}}}

\newcommand{\varCurriculumAtPeriod}{a^{\text{el}}}
\newcommand{\varCurriculumMandatoryAtPeriod}{a^{\text{mand}}}

\newcommand{\varCurMand}{\textnormal{vio}^{\text{mand}}}
\newcommand{\varCurThreePeriods}{\textnormal{vio}^{\text{3\_periods}}}
\newcommand{\varCurThreeMandPeriods}{\textnormal{vio}^{\text{3\_mand}}}
\newcommand{\varTutorium}{y^{\text{tut}}}
\newcommand{\varCourseOrderViolation}{\textnormal{vio}^{\text{order}}}
\newcommand{\varConsCourseSameRoomVio}{\textnormal{vio}^{\text{diff\_room}}}

\begin{document}

\title[A Multi-Objective Approach to Curriculum-Based Course Timetabling]{A Multi-Objective Approach to Curriculum-Based Course Timetabling with Continuity Across Semesters}

\author*[1]{\fnm{Florian} \sur{Meier}}\email{florian.meier@tum.de}

\author[2]{\fnm{Thomas} \sur{Stidsen}}\email{thst@dtu.dk}

\author[1]{\fnm{Clemens} \sur{Thielen}}\email{clemens.thielen@tum.de}

\author[3]{\fnm{Andreas} \sur{Wiese}}\email{andreas.wiese@tum.de}

\affil[1]{\orgdiv{Professorship of Optimization and Sustainable Decision Making, Campus Straubing for Biotechnology and Sustainability}, \orgname{Technical University of Munich},\orgaddress{
\city{Straubing}, \country{Germany}}}

\affil[2]{\orgdiv{Department of Technology, Management and Economics}, \orgname{Technical University of Denmark}, \city{Lyngby}, \country{Denmark}}

\affil[3]{\orgdiv{Department of Mathematics}, \orgname{Technical University of Munich}, \city{Munich}, \country{Germany}}

\abstract{We study a curriculum-based university course timetabling problem in which the preferences of two key stakeholder groups---lecturers and students---must be balanced while maintaining continuity across semesters in a weekly repeating timetable. While existing approaches typically rely on single-objective formulations or aggregate multiple objectives into a weighted sum, this can obscure the underlying trade-offs between conflicting stakeholder preferences. We therefore propose a multi-objective mixed-integer programming approach that explicitly separates lecturer and student objectives and incorporates timetable continuity by limiting the number of changes, called perturbations, in the time period assignments of selected courses relative to the corresponding semester of the previous academic year. To explore the resulting trade-offs, we develop a multi-objective solution approach based on the lexicographic $\varepsilon$-constraint method, enabling the computation of a representative set of solutions whose images, i.e., their vectors of objective values, cover different regions of the objective space.

The approach is evaluated on real-world instances from the Straubing Campus of the Technical University of Munich. The computational results reveal a clear and consistent trade-off between lecturers' and students' objectives across all instances. Moreover, the number of allowed perturbations is identified as a key decision parameter: relaxing this constraint significantly improves timetable quality for both stakeholder groups, although diminishing returns are observed beyond certain thresholds.

Overall, the proposed approach provides decision support by generating a diverse set of optimized timetables and enabling a transparent analysis of stakeholder trade-offs and continuity for the practical timetable planning process.
}

\keywords{University Course Timetabling, Mixed-Integer Programming, Multi-Objective Optimization, Stakeholder Trade-offs}

\maketitle

\section{Introduction}\label{UTT:sec:introduction}

Creating a timetable for an academic term is a fundamental, recurring task faced by universities worldwide~\citep{ceschia_educational_2023}. It leads to a university timetabling problem~\citep{bettinelli_overview_2015}, in which a given set of courses must be assigned to time periods and rooms while accounting for diverse requirements such as room capacities, lecturer availabilities, and structural restrictions arising from degree programs. 
In this paper, we study a curriculum-based course timetabling problem, a variant in which conflicts between courses are determined by predefined curricula. This differs from the related post-enrollment-based course timetabling problem, where conflicts are derived from the students' actual course enrollments. In our setting, a curriculum corresponds to a combination of a degree program and a semester of study, and the timetable should avoid overlaps among courses intended for the corresponding group of students.

In practice, constructing university timetables is often still a labor-intensive process requiring substantial manual work and iterative adjustments. Reducing this effort while improving timetable quality is an important operational task for universities. The constructed timetable directly shapes the students’ academic experience as well as the lecturers’ daily work routines. It determines whether the courses within each curriculum can be attended without timetable conflicts, and it influences how lecturers can coordinate teaching, research, and their other professional commitments.

University timetabling inherently involves multiple groups with differing and potentially conflicting interests~\citep{muhlenthaler_fairness_2016}. For example, several lecturers offering elective courses within the same curriculum may favor the same time periods during the week~\citep{babaei_applying_2018}, but this would result in overlaps within that curriculum that prevent students from attending (all of) these courses. Similarly, a timetable with breaks between lectures may benefit students by allowing recovery and preparation time, while lecturers who teach multiple of those courses may prefer to teach in consecutive time periods to free up larger uninterrupted time periods for research. Preferences regarding consecutively scheduled courses, distribution of courses across the week, preferred time periods, or simultaneously scheduled courses represent desirable timetable criteria rather than strict requirements, and they may not align between the stakeholder groups. It is not immediately clear to what extent these interests of the two main groups, students and lecturers, are compatible or conflicting, and how they should be handled in the timetabling process.

The relative importance of different timetable criteria of a stakeholder group can reasonably be assessed by the group members~\citep{babaei_applying_2018}. 
Although preferences may also vary within each stakeholder group~\citep{schimmelpfeng_application_2007}, lecturers can, for example, evaluate whether consecutive scheduling of their courses is more important overall than assigning them to preferred time periods, while students can judge the relative importance of avoiding late time periods versus minimizing overlaps of elective courses. It is, however, not straightforward how these criteria of the different stakeholder groups should be weighted against each other~\citep{vatandoost_ranking_2026}. Rather than being merely a modeling choice, it is a normative problem that directly influences the timetable outcomes. This raises the question of how such potential trade-offs should be treated in a principled manner.

An important practical requirement is continuity across semesters, i.e., that a course is scheduled in the same time period as in the last semester in which it was offered~\citep{lemos_disruptions_2021}. As most courses are offered on a yearly basis, these reference time periods usually originate from an earlier semester (typically the corresponding semester of the previous academic year, i.e., the previous summer semester or the previous winter semester) rather than the immediately preceding one.
While changes in room assignments are typically not critical, maintaining consistent time periods for selected courses is important, as academic activities, research collaborations, and recurring institutional commitments often extend beyond a single semester. Completely redesigning the timetables every year may disrupt established structures and coordination mechanisms. Therefore, the number of changes relative to the time period assignments in the corresponding semester of the previous academic year must be balanced carefully against improvements in timetable quality.

Given the considerations above, university timetabling should be treated as a multi-objective optimization problem~\citep{akkan_bi-criteria_2022}. It requires an optimization approach that allows different, potentially conflicting objectives to be considered simultaneously. Also, it should allow a transparent analysis of possible trade-offs, including the extent of permissible changes across different semesters. With such an optimization method, the timetable planners will be able to make better and more informed decisions regarding trade-offs. Also, they could more easily justify their decisions afterward to the different stakeholder groups.

\smallskip

The study presented in this paper is motivated by the real-world timetabling process at the TUM Campus Straubing for Biotechnology and Sustainability, which is a campus of the Technical University of Munich~(TUM). The underlying problem arises directly from the practical university timetable planning at the campus. There, a timetable that repeats on a weekly basis must be created for multiple degree programs with their specific structures. The problem setting was defined in close collaboration with the academic planners of the campus, and the data used in this paper stems from the real planning process. Timetables used in practice at the TUM Campus Straubing were computed with a single-objective version of the mixed-integer programming model presented in this paper, using a weighted sum of the different objectives. 
In addition, our work aims not only to address a practically relevant timetabling problem but also to provide methodological insights into how conflicting objectives can be modeled and analyzed in a real-world university environment.

The contributions of this paper can be summarized as follows. First, we propose a detailed mixed-integer programming (MIP) model for curriculum-based university timetabling that explicitly distinguishes between the two stakeholder groups, the students and the lecturers. Second, we develop a multi-objective solution approach that allows the computation of trade-off solutions that represent several possibilities to balance the interests of both groups. Third, we ensure continuity across semesters for selected courses by limiting the number of changes relative to their time period assignments in the corresponding semester of the previous academic year. Finally, we provide an empirical evaluation based on three real-world instances from the TUM Campus Straubing, thereby demonstrating both the practical relevance of the problem and the analytical value of the proposed multi-objective optimization approach.

The remainder of this paper is structured as follows. Section~\ref{UTT:sec:literature} provides an overview of the related literature. Section~\ref{UTT:sec:problem_description} presents a formal mathematical definition of the considered university timetabling problem, and Section~\ref{UTT:sec:MIP} introduces our corresponding MIP model. Section~\ref{UTT:sec:methodology} outlines the relevant concepts of multi-objective optimization (Section~\ref{UTT:subsec:multiobj}) and presents our multi-objective solution algorithm (Section~\ref{UTT:subsec:algorithm}). Section~\ref{UTT:sec:computational_results} describes the computational experiments and the real-world benchmark instances (Section~\ref{UTT:subsec:instances}) and reports the computational results (Section~\ref{UTT:sec:results}), including a final discussion and managerial insights. Finally, Section~\ref{UTT:sec:conclusion} concludes the paper. \section{Related Work}\label{UTT:sec:literature}

University course timetabling belongs to the broader field of educational timetabling, which is commonly divided into high-school timetabling, university examination timetabling, and university course timetabling (CTT)~\citep{ceschia_educational_2023}. Within this field, university course timetabling has received considerable attention in the discrete optimization and scheduling literature~\citep{chen_survey_2021}. Relevant problem variants include curriculum-based course timetabling (CB-CTT) and post-enrollment-based course timetabling. 
In the latter setting, the timetable is constructed after students have enrolled in specific courses, and the objective is to schedule courses in order to maximize the number of students who can attend all courses they have selected~\citep{jat_hybrid_2011}. By contrast, CB-CTT is performed \emph{before} enrollment, and conflicts are determined by predefined curricula published by the university~\citep{lach_curriculum_2012}. In this context, a curriculum typically corresponds to a combination of a degree program and a semester. Consequently, mandatory courses within the same curriculum must not be scheduled simultaneously, since they are intended for the same group of students, while overlaps among elective courses of the same curriculum are, in principle, tolerable but should be minimized. A broad range of solution methods has been proposed, including heuristics and metaheuristics such as tabu search, simulated annealing, and genetic algorithms~\citep{chiarandini_effective_2006, bellio_design_2012, dunke_matheuristic_2023}, as well as exact MIP-based approaches~\citep{schimmelpfeng_application_2007, lach_curriculum_2012}.

In this paper, we focus on a CB-CTT problem, as this is the problem variant that arises in the real-world timetabling process at the TUM Campus Straubing. CB-CTT problems have been studied in several variants. These include, for example, formulations that assess timetable quality at the level of individual curricula rather than optimizing the overall timetable quality across all curricula~\citep{muhlenthaler_fairness_2016}, as well as extensions such as student sectioning, where students must be assigned to specific sessions of a course if that course is offered multiple times, as considered in the International Timetabling Competition~(ITC) in 2019~\citep{muller_real-world_2025}. Among the existing variants, the problem definition of the ITC 2007~\citep{mccollum_setting_2010} is one of the most widely studied variants, and its instances have served as benchmarks throughout the literature~\citep{lindahl_quality_2019, gulcu_robust_2020, lach_curriculum_2012, muhlenthaler_fairness_2016}. However, this formulation was deliberately simplified with respect to real-world university timetabling for research and competition purposes~\citep{mccollum_setting_2010}. At the ITC~2007, solution quality was evaluated by penalizing violations of four soft constraint categories—room capacity, minimum working days, curriculum compactness, and room stability—assigning a fixed weight to each of these criteria.
The different solution approaches for CB-CTT problems in the literature often rely on a single aggregated objective~\citep{lach_curriculum_2012}. As a result, they are less suited to settings where it is unclear which weight should be given to which criterion in a timetable, and where the trade-offs between such criteria should be analyzed explicitly.

An important practical aspect in university timetabling is continuity across semesters. In this context, it is important to distinguish between robustness and stability. Robust timetabling approaches aim to construct timetables that remain feasible under potential disruptions, such as short-notice changes in resource availability, by minimizing the number of assignments that need to be modified when disruptions occur. For example, \citet{akkan_bi-criteria_2022} and \citet{gulcu_robust_2020} develop approaches that explicitly account for such uncertainty during timetable construction.
In contrast, stability and minimal perturbation approaches focus on adapting an existing timetable after disruptions have occurred. In the minimal perturbation setting, the objective is to compute a new timetable that deviates as little as possible from a previously established one, typically by minimizing the number of changes, as studied by \citet{phillips_integer_2017} and \citet{burke_minimal_2005}. Stability-oriented approaches extend this idea by additionally considering timetable quality. For instance, \citet{perzina_solving_2007} and \citet{lemos_disruptions_2021, hebrard_minimal_2020} incorporate solution quality either lexicographically or via weighted sums, while \citet{lindahl_quality_2019} explicitly analyze the trade-off between timetable quality and the number of perturbations. In contrast to these approaches focusing on reacting to disruptions, we consider a different setting. We address the construction of a timetable for a new semester, where continuity is defined with respect to the most recent semester in which each course was offered. In particular, courses should retain their previously assigned time periods whenever this is desired, while changes in room assignments are considered irrelevant. Thus, we incorporate continuity as a planning objective that selectively preserves time period assignments across semesters.

Some multi-objective approaches for university course timetabling problems exist in the literature~\citep{Datta2007, davison_modelling_2025} and some robustness- and stability-oriented approaches discussed above are multi-objective optimization methods, as they balance timetable quality with robustness or the number of perturbations.
However, explicit distinctions between stakeholder groups such as students and lecturers remain rare in the literature, and the corresponding trade-offs are typically not analyzed directly. Although \citet{dunke_matheuristic_2023} incorporate different stakeholders' perspectives in a multi-level, multi-criteria timetabling framework, they do not distinguish between lecturers' and students' objectives and neither consider continuity across semesters. 
Hence, the explicit multi-objective modeling and analysis of students’ and lecturers’ objectives, combined with continuity across semesters, remains largely unaddressed.

To summarize, the existing literature addresses university course timetabling from a variety of perspectives, including multi-objective optimization. However, the trade-offs between lecturers and students are not explicitly captured in the existing literature. In particular, these trade-offs have not been analyzed in combination with continuity across semesters. To the best of our knowledge, continuity across semesters has not been incorporated into the construction of new timetables in the way considered in this paper, where continuity across semesters is modeled for selected courses by limiting the number of deviations from their previously assigned time periods, while room assignments may change. \section{Problem Description}\label{UTT:sec:problem_description}

The multi-objective CB-CTT problem we consider concerns generating a timetable for a given set~$\setCourses$ of university courses. This problem focuses on creating a timetable for a specific semester, where the timetable repeats on a weekly basis. A given set~$\setLecturers$ of lecturers and a given set~$\setRooms$ of rooms are available. Each course~$c\in\setCourses$ is associated with, among other attributes, a set~$\setLecturersOfCourse(c)\subseteq\setLecturers$ of lecturers teaching it, the expected number~$\parSeatsForCourse(c)$ of attending students, which determines the number of seats required in any assigned room, and a specified number of sessions denoted by $\parTimesPerWeek(c)\in\mathbb{N}$, indicating the number of times it must be scheduled within a week. Each course~$c\in\setCourses$ must be allocated to $\parTimesPerWeek(c)$ 2-hour time slots, hereafter called periods. Each working day~$d \in \setDays$ is divided into six periods between 8~a.m. and 8~p.m., and the set~$\setPeriodsAtDay(d)$ denotes the periods on day~$d$. Altogether, the set~$\setPeriods = \bigcup_{d\in \setDays} \setPeriodsAtDay(d)$ of all periods contains 30 elements. Moreover, we denote by~$\setSubsequentPeriods\subseteq\setPeriods$ the set of periods that have a subsequent period on the same day, i.e., all periods except the last period of each day.
The goal is to determine a timetable in which each course~$c \in \setCourses$ is assigned to~$\parTimesPerWeek(c)$ period-room pairs, where each pair specifies that the course takes place in the corresponding room during the corresponding period. 

The set~$\setCourses$ of courses is partitioned into modules~$M\in\setModules$, i.e., each course~$c\in\setCourses$ belongs to exactly one module. Each module consists of closely related courses that jointly cover a common topic (e.g., a lecture and its associated tutorials) and therefore must be scheduled in a coordinated manner. Therefore, we distinguish between three different course types, which induce a second partition of the course set into the set~$\setLectures\subseteq\setCourses$ of lectures, the set~$\setExercises\subseteq\setCourses$ of (plenary) exercise classes, and the set~$\setTutorials\subseteq\setCourses$ of tutorials. Exercise classes are plenary formats in which all students attend the same class. In contrast, a tutorial  
may consist of multiple sessions scheduled either in parallel or at different time periods, covering the same material. Students are expected to attend exactly one of these sessions.
Typically, in a tutorial students work on exercises in small groups and they can clarify questions. For each module, at most one tutorial course is offered.

The timetable must satisfy a set of constraints related to room assignments. These include that the seating capacity of a room~$r$, denoted by~$\parSeatsInRoom(r)$, must be at least the number~$\parSeatsForCourse(c)$ of seats required for any course~$c$ that~$r$ is assigned to. Other requirements may include the availability of computers or power outlets for students. For each course~$c\in\setCourses$, the set~$\setRequiredRooms(c)\subseteq\setRooms$ contains the rooms that satisfy all requirements of~$c$. Conversely, for each room~$r\in\setRooms$, the set~$\setCompatibleCourses(r)\subseteq\setCourses$ contains all courses whose requirements are satisfied by~$r$. 
Further, we introduce the set~$\setRoomConflicts$ of room conflicts, which consists of two-element subsets of~$\setRooms$, where $\{r_1,r_2\}\in\setRoomConflicts$ indicates that rooms~$r_1$ and~$r_2$ cannot be used simultaneously. 
For example, some rooms can be subdivided using partition walls, allowing a larger room to be split into two smaller rooms. In that case, each smaller room~$r_1\in\setRooms$ resulting from the subdivision can only be used when the larger room~$r_2\in\setRooms$ is not being used as a whole.
Because assembling and removing the partition walls is a time-consuming process, a decision must be made each day~$d\in\setDays$ on whether to use the room as two smaller rooms or keep it as a large one. We refer to a room as \emph{active} on a given day if it is enabled for assignment throughout that day. The set~$\setRoomDaySetting\subseteq\setRooms$ contains the rooms for which a daily activation decision must be made: on each day~$d$, such a room is either active in all periods~$p\in\setPeriodsAtDay(d)$ or inactive throughout the day. 
Additionally, period-specific restrictions may apply even to rooms that are active on a given day. In particular, some rooms may be unavailable in certain periods because they are rented for external events. Thus, for each room~$r\in\setRooms$, the set~$\setExcludedPeriodsOfRoom(r)\subseteq\setPeriods$ contains all periods in which no course may be scheduled in that room.

For both lecturers and students, it is desirable that lectures take place in the same time periods as in the previous semester in which the corresponding course was offered, e.g., for planning of future courses and external commitments. Therefore, for each course~$c \in \setCourses$, the set~$\setPeriodsOfLastYear(c)\subseteq\setPeriods$ contains the periods in which the course was scheduled in the previous year if continuity across semesters is desired for that course, and is empty otherwise. We call an assignment of a course~$c \in \setCourses$ with~$\setPeriodsOfLastYear(c)\neq\emptyset$ to a period outside~$\setPeriodsOfLastYear(c)$ a \emph{perturbation}. The parameter~$\parPerturbationsBound$ specifies an upper bound on the number of perturbations, i.e., on the number of such assignments to periods outside the respective previous periods.

The main goal of the timetabling problem is to accommodate the preferences of different stakeholder groups, namely lecturers and students. For instance, lecturers may favor certain periods for their courses, while students may wish to avoid classes on Friday evenings or overlaps between different courses that are relevant for their curriculum. However, since the timetable comprises only 30 periods and most curricula include more than 30 courses, it is generally impossible to fully satisfy all preferences. The constraints and objectives associated with both stakeholder groups are defined below.

\subsection{Lecturers}

The availability and the preferences of the lecturers play an important role in the scheduling process.
For example, the scheduling must account for periods during which lecturers are unavailable due to other commitments. Therefore, for each lecturer $l\in\setLecturers$, we denote by $\setCoursesOfLecturer(l)\subseteq\setCourses$ the set of courses taught by~$l$, and we introduce the set~$\setExcludedPeriods(l)\subsetneq\setPeriods$ of excluded periods in which these courses cannot be scheduled.
Moreover, some courses $c\in\setCourses$ must be scheduled in predefined periods~$\setfixPeriods(c)\subseteq\setPeriods$ because, for example, they are in cooperation with external lecturers.

To incorporate preferred periods in the computed timetable, we define for each course~$c \in \setCourses$ the set~$\setPrefPeriodsOfCourse(c) \subseteq \setPeriods$ of preferred periods, which comprises all periods favored by the lecturers teaching course~$c$.
Some lecturers prefer their courses to be scheduled in consecutive periods within a teaching day and in the same room; these lecturers are collected in the set~$\setLecturersCompactPeriods\subseteq\setLecturers$. Lecturers who instead prefer their courses to be scheduled in non-consecutive periods, with a break between teaching periods, are included in the set~$\setLecturersDistPeriods\subseteq\setLecturers$. In addition, lecturers who indicate a preference for teaching on consecutive days are assigned to the set~$\setLecturersCompactWeek\subseteq\setLecturers$, whereas those who prefer teaching days to be separated by at least one free day are included in~$\setLecturersDistWeek\subseteq\setLecturers$, with~$\setLecturersCompactWeek\cap\setLecturersDistWeek=\emptyset$. A lecturer may also belong to none of these sets if they are indifferent and express no specific preference.
Each lecturer~$l\in\setLecturers$ may also specify a preferred number~$\parNumberOfDays(l)$ of teaching days over which their courses should be distributed equally. Lecturers who have indicated such a preference are included in the set~$\setLecturersPrefDaysPerWeek\subseteq\setLecturers$.

To ensure fairness, a specific number can be set as a limit on how often lecturers are assigned to teach at certain times during the week. For example, a lecturer may be scheduled to teach in the 6~p.m. to 8~p.m. period at most once per week. Accordingly, a subset~$G\subseteq\setPeriods$ of periods is called a \emph{period group} if each lecturer may be assigned to teach in at most~$\parMaxPeriodSet(G)$ periods in~$G$. The set of all period groups is denoted by~$\setPeriodGroups$.

\subsection{Students}
In addition to lecturers, we consider students as the second stakeholder group.  Each student is enrolled in a specific degree program and is in a particular semester of study. To account for the objectives of these students, we model each combination of degree program and semester of study (e.g., Bachelor Bioeconomy, third semester) as a distinct curriculum, and we denote by~$\setCurricula$ the set of all curricula. The set of courses associated with curriculum~$\kappa\in\setCurricula$ is denoted by~$\setOfCoursesInCurr(\kappa)\subseteq\setCourses$. These courses should be scheduled in a manner that enables the corresponding student cohort to attend them and results in a well-structured timetable.

For each curriculum~$\kappa \in \setCurricula$, the associated courses are classified as either elective or mandatory, with the corresponding sets denoted by~$\setElectiveCourses(\kappa)$ and~$\setMandCourses(\kappa)$, respectively. Each elective course~$c\in\setElectiveCourses(\kappa)$ may be chosen by students of curriculum~$\kappa$ and is assigned a weight~$\parCourseCurrWeight(c,\kappa)$, where higher weights indicate greater importance within the curriculum and reflect that more students are expected to choose the course. Mandatory courses, in contrast, must be attended by all students in curriculum~$\kappa$.
A key requirement of the timetable is to avoid overlaps between mandatory courses within the same curriculum; that is, mandatory courses must not be scheduled in the same period. In addition, the timetable should minimize the number of overlaps between mandatory and elective courses, as well as among elective courses themselves.
In the case of tutorials, students are only required to be able to attend one session. Therefore, for each tutorial~$c \in \setTutorials$ and each curriculum~$\kappa \in \setCurricula$ with~$c \in \setMandCourses(\kappa)$, it is sufficient that at least one session of~$c$ does not overlap with a mandatory course of curriculum~$\kappa$; the remaining sessions may overlap. 

Ideally, breaks should be ensured between the courses in a curriculum to give students time for moving between rooms, short recovery, and preparation. This may naturally contradict the objective of lecturers who prefer consecutive teaching periods, since a single lecturer often teaches multiple courses within the same curriculum and may therefore wish to have those courses scheduled consecutively. Breaks are therefore considered in the students' objective.

Further, courses within the same module may be subject to additional structural requirements. Some courses should, for example, be scheduled in a specific order within the week. In particular, the exercise class of a module should usually take place after the corresponding lecture within the same week, so that it can build on the material covered there. To model this, for each module~$M\in\setModules$, we are given a relation~$\prec$ on the courses in~$M$, where~$c_1 \prec c_2$ for~$c_1,c_2\in M$ indicates that the first session of~$c_1$ should be scheduled before all sessions of~$c_2$. In this context, it may also be required or desirable that two courses be scheduled on different days. 
To capture this, we introduce the set~$\setCourseConflict$ of hard course conflicts and the set~$\setCourseTuples$ of soft course conflicts. A subset $\{c_1,c_2\}\in\setCourseConflict$ indicates that courses~$c_1$ and~$c_2$ must not be scheduled on the same day, whereas $\{c_1,c_2\}\in\setCourseTuples$ indicates that scheduling them on the same day is possible but undesired.
Moreover, for any course $c \in \setCourses$, a follow-up course~$\parFollowingCourse(c)$ may be specified that has to be scheduled immediately after~$c$ in the next time slot on the same day and in the same room. This may be desirable, for example, when an exercise class is intended to directly build on the content of a lecture, allowing students to apply the material immediately in exercises. 

Additionally, some periods, such as Friday evenings from 6~p.m. to 8~p.m., are generally unpopular with students and should therefore be avoided whenever possible. Accordingly, assigning a course~$c$ to a period~$p\in\setPeriods$ incurs a penalty proportional to the period-specific penalty parameter~$\parPeriodPenalty(p)$ and the number $\parSeatsForCourse(c)$ of students expected to attend the course. Moreover, to preserve rooms for independent student study, certain rooms should preferably remain unassigned, although they may be used if necessary. To discourage such use, assigning a course to room~$r$ incurs a per-period penalty given by~$\parRoomPenalty(r)$. \section{Mathematical Model}\label{UTT:sec:MIP}

This section presents the mathematical formulation of the objective functions and constraints of the MIP model used to solve the CB-CTT problem described in the previous section. We first define all sets and parameters. Sets with superscripts are always subsets of the set denoted by the same symbol without superscript, for example, $\setLectures\subseteq\setCourses$.

\subsection{Sets and parameters}

In the following, all sets and parameters introduced in Section~\ref{UTT:sec:problem_description} and used in the MIP formulation are listed for convenience. Sets are denoted by calligraphic capital letters, whereas parameters are denoted by lowercase symbols or descriptive textual notation.

\medskip

\begin{xltabular}{\textwidth}{lX}
$\setCourses$ &set of courses (index~$c$)\\
$\setLecturers$ &set of lecturers (index~$l$)\\
$\setRooms$ &set of rooms (index~$r$)\\
$\setLecturersOfCourse(c)$ &set of lecturers teaching course~$c\in\setCourses$\\
$\parSeatsForCourse(c)\in\mathbb{N}$ & number of seats required for course~$c\in\setCourses$\\
$\parTimesPerWeek(c)\in\mathbb{N}$ & number of sessions of course~$c\in\setCourses$ per week\\
$\setDays=\{1,\dots,5\}$ &set of five working days (index~$d$)\\
$\setPeriods=\{1,\dots,30\}$ & set of periods over the week, for five days with six periods each day (index~$p$)\\
$\setPeriodsAtDay(d)$ &set of periods on day~$d\in\setDays$\\
$\setSubsequentPeriods$ &set of periods with a subsequent period on the same day, i.e., excluding each day’s last period\\
$\setModules$ &set of modules, where each module~$M\in\setModules$ is a subset~$M\subseteq\setCourses$ of courses (index~$M$)\\
$\setLectures$ &set of lectures\\
$\setExercises$ &set of (plenary) exercise classes\\
$\setTutorials$ &set of tutorials\\
$\parSeatsInRoom(r)\in\mathbb{N}$ &number of seats in room~$r\in\setRooms$\\
$\setRequiredRooms(c)$ &set of rooms that meet the requirements of course~$c\in\setCourses$\\
$\setCompatibleCourses(r)$ &set of all courses whose requirements are satisfied by room~$r$\\
$\setRoomConflicts$ &set of room conflicts, where $\{r_1, r_2\} \in \setRoomConflicts$ indicates that rooms~$r_1$ and~$r_2$ cannot be used simultaneously\\
$\setRoomDaySetting$ &set of rooms for which a daily activation decision must be made, such that each room is either active in all periods~$p\in\setPeriodsAtDay(d)$ of a day~$d\in\setDays$ or in none of them\\
$\setExcludedPeriodsOfRoom(r)$ &set of periods during which room~$r\in\setRooms$ is unavailable\\
$\setCoursesOfLecturer(l)$ &set of courses taught by lecturer~$l\in\setLecturers$\\
$\setExcludedPeriods(l)$ &set of excluded periods of lecturer~$l\in\setLecturers$ during which they cannot teach\\
$\setfixPeriods(c)$ &set of predefined periods during which course~$c$ must take place\\
$\setPrefPeriodsOfCourse(c)$ &set of preferred periods of course~$c\in\setCourses$\\
$\setLecturersCompactPeriods$ &set of lecturers who prefer their courses to be scheduled in consecutive periods within a teaching day\\
$\setLecturersDistPeriods$ &set of lecturers who prefer their courses to be scheduled in non-consecutive periods within each teaching day\\
$\setLecturersCompactWeek$ &set of lecturers who prefer teaching on consecutive days\\
$\setLecturersDistWeek$ &set of lecturers who prefer teaching days to be separated by at least one free day\\
$\setLecturersPrefDaysPerWeek$ &set of lecturers with a preferred number of teaching days per week\\
$\parNumberOfDays(l)$ &preferred number of teaching days per week of lecturer~$l\in\setLecturersPrefDaysPerWeek$\\
$\setPeriodGroups$ &set of period groups, where each period group~$G \in \setPeriodGroups$ is a subset~$G\subseteq\setPeriods$ of periods\\
$\parMaxPeriodSet(G)$ &maximum number of periods in period group~$G\in\setPeriodGroups$ in which each lecturer may teach\\
$\setPeriodsOfLastYear(c)$ &set of periods in which course~$c\in\setCourses$ was scheduled in the previous year if continuity across semesters is desired for this course, otherwise $\setPeriodsOfLastYear(c)=\emptyset$\\
$\setCurricula$ &set of curricula (index~$\kappa$)\\
$\setOfCoursesInCurr(\kappa)$ &set of courses associated with curriculum~$\kappa \in \setCurricula$\\
$\setElectiveCourses(\kappa)$ &set of elective courses in curriculum~$\kappa \in \setCurricula$\\
$\setMandCourses(\kappa)$ &set of mandatory courses in curriculum~$\kappa \in \setCurricula$\\
$\parCourseCurrWeight(c, \kappa)\in\mathbb{N}$ & relative importance of elective course~$c\in\setElectiveCourses(\kappa)$ within curriculum~$\kappa\in\setCurricula$\\
$\setCourseConflict$ &set of hard course conflicts, where $\{c_1,c_2\}\in\setCourseConflict$ indicates that courses~$c_1$ and~$c_2$ must not be scheduled on the same day\\
$\setCourseTuples$ &set of soft course conflicts, where $\{c_1,c_2\}\in\setCourseTuples$ indicates that courses~$c_1$ and~$c_2$ should preferably not be scheduled on the same day\\
$\parFollowingCourse(c)$ & follow-up course that has to be scheduled immediately after course~$c\in\setCourses$ in the next time slot on the same day and in the same room. This parameter is set to \texttt{None} if no follow-up course is specified for course~$c$\\
$\parPeriodPenalty(p)\in\mathbb{Q}_{\geq 0}$ &penalty from the students' perspective for using period~$p\in\setPeriods$\\
$\parRoomPenalty(r)\in\mathbb{Q}_{\geq 0}$ &per-period penalty for assigning a course to room~$r\in\setRooms$\\
$\parLecturersBound\in\mathbb{Q}\cup\{+\infty\}$ &upper bound on the lecturers' objective\\
$\parStudentsBound\in\mathbb{Q}\cup\{+\infty\}$ &upper bound on the students' objective\\
$\parPerturbationsBound\in\mathbb{N}\cup\{+\infty\}$ &upper bound on the number of perturbations
\end{xltabular}

\subsection{Variables}
We use three groups of variables, denoted by~$a$, $y$, and~$\textnormal{vio}$ with corresponding indices. 
The assignment variables, denoted by~$a$, define the overall timetable structure, while the $y$-variables serve as auxiliary variables. The variables labeled $\textnormal{vio}$ are used to penalize soft constraint violations. In the following, we denote the vector comprising all variables of our MIP model by~$x$.

\medskip

\begin{xltabular}{\textwidth}{lX}
   $\varAssignment_{c,p,r}$ &binary variable equal to~1 if course~$c\in\setCourses$ is assigned to period~$p\in\setPeriods$ and room~$r\in\setRequiredRooms(c)$; 0~otherwise\\
   $\varLecAtPeriod_{l, p}$ &binary variable equal to~1 if lecturer~$l\in\setLecturers$ teaches in period~$p\in\setPeriods$; 0~otherwise\\
   $\varLecAtDay_{l, d}$ &binary variable equal to~1 if lecturer~$l\in\setLecturers$ teaches on day~$d\in\setDays$; 0~otherwise\\
   $\varCurriculumAtPeriod_{\kappa, p}$ &binary variable equal to~1 if an elective lecture or exercise class of curriculum~$\kappa\in\setCurricula$ is assigned to period~$p\in\setPeriods$; 0~otherwise\\
   $\varCurriculumMandatoryAtPeriod_{\kappa, p}$ &binary variable equal to~1 if a mandatory lecture or exercise class of curriculum~$\kappa\in\setCurricula$ is assigned to period~$p\in\setPeriods$; 0~otherwise\\    
   $\varRoomAvailability_{r, p}$ &binary variable equal to~1 if room~$r\in\setRoomDaySetting$ is active in period~$p\in\setPeriods$; 0 otherwise\\
   $\varTutorium_{\kappa, c, p}$ & binary variable equal to~1 if period~$p\in\setPeriods$ is selected as the period in which students of curriculum~$\kappa\in\setCurricula$ can attend one session of the mandatory tutorial~$c\in\setMandCourses(\kappa)\cap\setTutorials$; 0~otherwise\\
   $\varLecConsecutivePeriods_{l, p}$ &binary variable equal to~1 if lecturer~$l\in\setLecturersDistPeriods$ teaches in the two consecutive periods~$p\in\setSubsequentPeriods$ and~$p+1$ in one day; 0~otherwise\\
   $\varLecDisjointPeriods_{l, p}$ & binary variable equal to~1 if lecturer~$l\in\setLecturersCompactPeriods$ teaches in period~$p\in \bigcup_{d\in\setDays}\{\min(\setPeriodsAtDay(d)), \dots , \max(\setPeriodsAtDay(d)) -2 \}$ and in a later period on the same day, with at least one free period in between; 0~otherwise\\
   $\varLecAtDayPlusOne_{l, d}$ &binary variable equal to~1 if lecturer~$l\in\setLecturersDistWeek$ teaches on day~$d\in\{1,\dots,\lvert\setDays \rvert - 1\}$ and on day~$d+1$; 0~otherwise\\
   $\varLecAtDayPlusTwo_{l, d}$ & binary variable equal to~1 if lecturer~$l\in\setLecturersCompactWeek$ teaches on day~$d\in{1,\dots,\lvert\setDays \rvert - 2}$ and again on day~$d+2$, with no teaching on the intervening day~$d+1$; 0~otherwise\\
   $\varLecAtDayPlusThree_{l, d}$ & binary variable equal to~1 if lecturer~$l\in\setLecturersCompactWeek$ teaches on day~$d\in{1,\dots,\lvert\setDays \rvert - 3}$ and again on day~$d+3$, with no teaching on the two intervening days~$d+1$ and~$d+2$; 0~otherwise\\
   $\varLecAtDayPlusFour_{l, d}$ & binary variable equal to~1 if lecturer~$l\in\setLecturersCompactWeek$ teaches on day~$d\in{1,\dots,\lvert\setDays \rvert - 4}$ and again on day~$d+4$, with no teaching on the three intervening days~$d+1$, $d+2,$ and~$d+3$; 0~otherwise\\
   $\varLecDevDaysPos_{l}$ & nonnegative integer variable representing the difference between the actual number of teaching days of lecturer~$l\in\setLecturersPrefDaysPerWeek$ and the desired value~$\parNumberOfDays(l)$, if the actual number is higher; 0~otherwise\\
   $\varLecDevDaysNeg_{l}$ & nonnegative integer variable representing the difference between the desired value~$\parNumberOfDays(l)$ and the actual number of teaching days of lecturer~$l\in\setLecturersPrefDaysPerWeek$, if the desired value is higher; 0~otherwise\\
   $\varLecDeviationCoursesPerDay_{l, d}$ &nonnegative integer variable representing the absolute deviation of the number of courses assigned to lecturer~$l\in\setLecturersPrefDaysPerWeek$ on teaching day~$d\in\setDays$ from the desired average number of courses for each teaching day\\
   $\varConsCourseSameRoomVio_{l, p}$ &binary variable equal to~1 if lecturer~$l\in\setLecturers$ teaches in the two consecutive periods $p\in\setSubsequentPeriods$ and~$p+1$, but in different rooms; 0~otherwise\\
   $\varCourseDistVio_{c_1, c_2, d}$ & binary variable equal to~1 if courses~$c_1,c_2\in\setCourses$ with~$\{c_1,c_2\}\in\setCourseTuples$ are both scheduled on day~$d$, even though they should preferably not be scheduled on the same day; 0~otherwise\\
   $\varCourseOrderViolation_{c_1, c_2}$ & binary variable equal to~1 if courses~$c_1,c_2\in\setCourses$ satisfy~$c_1\prec c_2$ but are scheduled in the opposite order within the week; 0~otherwise\\
   $\varCurCourseVio_{\kappa, c, p}$ &binary variable equal to~1 if an elective lecture or exercise class~$c\in(\setLectures \cup \setExercises) \cap \setElectiveCourses(\kappa)$ of curriculum~$\kappa\in\setCurricula$ overlaps with another course of curriculum~$\kappa$ in period~$p\in\setPeriods$; 0~otherwise\\
   $\varCurMand_{\kappa, p}$ &binary variable equal to~1 if two mandatory lectures or exercise classes of curriculum~$\kappa\in\setCurricula$ are assigned to periods $p\in\setSubsequentPeriods$ and $p+1$; 0~otherwise\\
   $\varCurThreePeriods_{\kappa, p}$ & binary variable equal to~1 if, in curriculum~$\kappa\in\setCurricula$, two mandatory and one elective lecture or exercise class are scheduled in period~$p\in\setPeriods$ and the two subsequent periods on the same day; 0~otherwise\\
   $\varCurThreeMandPeriods_{\kappa, p}$ & binary variable equal to~1 if, in curriculum~$\kappa\in\setCurricula$, three mandatory lectures or exercise classes are scheduled in period~$p\in\setPeriods$ and the two subsequent periods on the same day; 0~otherwise\\
\end{xltabular}

\subsection{Objectives}
We consider three objectives in the CB-CTT problem, namely the lecturers' objective~$\parLecturersObj$, the students' objective~$\parStudentsObj$, and the perturbations objective~$\parPerturbationsObj$. The lecturers' and students' objectives are optimized simultaneously in a bi-objective optimization approach in order to identify trade-offs between the two stakeholder groups. The perturbations objective counts the number of perturbations, i.e., the number of session assignments of courses~$c \in \setCourses$ with~$\setPeriodsOfLastYear(c)\neq\emptyset$ to periods outside~$\setPeriodsOfLastYear(c)$. It is first minimized to determine the minimum number of unavoidable changes and subsequently bounded from above by different values to limit the number of permitted changes. The lecturers' and students' objectives are each defined as weighted sums of several subobjectives, where the weights~$\alpha_i$ used in the following are described in Appendix~\ref{UTT:app:survey}.

\subsubsection{Lecturers' Objective}\label{UTT:lecturers_obj}

The lecturers' objective~$\parLecturersObj$ is defined as $\parLecturersObj(x) = \sum_{i=1}^{10} \alpha_i z_i(x)$, where the subobjectives~$z_1$ to~$z_{10}$, defined below, capture the lecturers' preferences and are all to be minimized.

\medskip

\noindent

\noindent
Assignment of courses to periods that are not preferred by the lecturers teaching the course, taking into account the number of lecturers involved in the course:
\begin{align}
z_1(x) = \sum_{c \in C : \setPrefPeriodsOfCourse(c)\neq\emptyset} \sum_{p\in\setPeriods\setminus\setPrefPeriodsOfCourse(c)} \sum_{r \in \setRequiredRooms(c)} \left( 1 + \frac{\lvert\setLecturersOfCourse(c)\rvert - 1}{3} \right) \cdot\varAssignment_{c,p,r}
\end{align}

\noindent
Consecutive teaching periods of a lecturer within a day although the lecturer prefers non-consecutive teaching periods:
\begin{align}
z_2(x)=\sum_{l \in \setLecturersDistPeriods} \sum_{p \in \setPeriods} \varLecConsecutivePeriods_{l, p}  
\end{align}

\noindent
Non-consecutive teaching periods of a lecturer within a day although the lecturer prefers consecutive teaching periods:
\begin{align}
z_3(x)=\sum_{l \in \setLecturersCompactPeriods} \sum_{p \in \setPeriods} \varLecDisjointPeriods_{l, p}  
\end{align}

\noindent
Consecutive teaching days of a lecturer although the lecturer prefers teaching days to be separated by at least one free day:
\begin{align}
z_4(x)=\sum_{l \in \setLecturersDistWeek} \sum_{d \in \{ 1, \dots, \lvert\setDays \rvert - 1 \}} \varLecAtDayPlusOne_{l, d}  
\end{align}

\noindent
Teaching on day~$d\in \{ 1, \dots, \lvert\setDays \rvert - 2 \}$ and again on day~$d+2$, with no teaching on the intervening day, although the lecturer prefers consecutive teaching days:
\begin{align}
z_5(x)=\sum_{l \in \setLecturersCompactWeek} \sum_{d \in \{ 1, \dots, \lvert\setDays \rvert - 2 \}} \varLecAtDayPlusTwo_{l, d}  
\end{align}

\noindent
Teaching on day~$d\in \{ 1, \dots, \lvert\setDays \rvert - 3 \}$ and again on day~$d+3$, with no teaching on the two intervening days, although the lecturer prefers consecutive teaching days:
\begin{align}
z_6(x)=\sum_{l \in \setLecturersCompactWeek} \sum_{d \in \{ 1, \dots, \lvert\setDays \rvert - 3 \}} \varLecAtDayPlusThree_{l, d} 
\end{align}

\noindent
Teaching on day~$d\in \{ 1, \dots, \lvert\setDays \rvert - 4 \}$ and again on day~$d+4$, with no teaching on the three intervening days, although the lecturer prefers consecutive teaching days:
\begin{align}
z_7(x)=\sum_{l \in \setLecturersCompactWeek} \sum_{d \in \{ 1, \dots, \lvert\setDays \rvert - 4 \}} \varLecAtDayPlusFour_{l, d} 
\end{align}

\noindent
Deviations from the desired number of teaching days of a lecturer:
\begin{align}
z_8(x)=\sum_{l \in \setLecturersPrefDaysPerWeek} \varLecDevDaysPos_{l} + \varLecDevDaysNeg_{l} 
\end{align}

\noindent
Deviations of the number of courses on a teaching day of a lecturer from the desired average number of courses per teaching day:
\begin{align}
z_{9}(x)=\sum_{l \in \setLecturersPrefDaysPerWeek} \sum_{d \in \setDays} \varLecDeviationCoursesPerDay_{l, d}  
\end{align}

\noindent
Consecutive courses of a lecturer that are not follow-up courses but are assigned to different rooms:
\begin{align}
z_{10}(x)=\sum_{l \in \setLecturers} \sum_{p \in \setSubsequentPeriods} \varConsCourseSameRoomVio_{l, p}
\end{align}

\subsubsection{Students' Objective}\label{UTT:students_obj}
The students' objective~$\parStudentsObj$ is defined as $\parStudentsObj(x) = \sum_{i=11}^{19} \alpha_i\cdot z_i(x)$, where the subobjectives~$z_{11}$ to~$z_{19}$, defined below, capture the students' preferences and are all to be minimized.

\medskip

\noindent
Assignment of courses that should preferably be scheduled on different days but are both scheduled on the same day:
\begin{align}
z_{11}(x)= \sum_{\{c_1, c_2\}\in\setCourseTuples} \sum_{d \in D} \varCourseDistVio_{c_1, c_2, d} 
\end{align}

\noindent
Courses within the same module that satisfy~$c_1 \prec c_2$ but are scheduled in the opposite order within the week:
\begin{align}
z_{12}(x)=\sum_{M\in\setModules} \quad \sum_{c_1, c_2 \in M: c_2 \prec c_1} \varCourseOrderViolation_{c_1, c_2}
\end{align}

\noindent
Consecutive scheduling of two mandatory lectures or exercise classes within a curriculum:
\begin{align}
z_{13}(x)=\sum_{\kappa \in \setCurricula} \sum_{p \in \setSubsequentPeriods} \varCurMand_{\kappa, p}  
\end{align}

\noindent
Consecutive scheduling of two mandatory and one elective lecture or exercise class within a curriculum:
\begin{align}
z_{14}(x)=\sum_{\kappa \in \setCurricula} \sum_{p \in \setPeriods} \varCurThreePeriods_{\kappa, p}  
\end{align}

\noindent
Consecutive scheduling of three mandatory lectures or exercise classes within a curriculum:
\begin{align}
z_{15}(x)=\sum_{\kappa \in \setCurricula} \sum_{p \in \setPeriods} \varCurThreeMandPeriods_{\kappa, p}  
\end{align}

\noindent
Overlap of an elective course $c\in\setElectiveCourses(\kappa)$ with any other course in curriculum~$\kappa\in\setCurricula$, where the courses considered are lectures or exercise classes, weighted by the importance~$\parCourseCurrWeight(c, \kappa)$ of course~$c$ within curriculum~$\kappa$:
\begin{align}
z_{16}(x)=\sum_{\kappa \in \setCurricula} \sum_{p \in \setPeriods} \sum_{c \in \setElectiveCourses(\kappa) \cap (\setLectures \cup \setExercises)} \parCourseCurrWeight(c, \kappa) \cdot \varCurCourseVio_{\kappa, c, p}
\end{align}

\noindent
Differences between the number of seats required for a course and the seating capacity of the assigned room. This should avoid, for example, using large rooms for small tutorials:
\begin{align}
z_{17}(x)=\sum_{c \in \setTutorials} \sum_{r \in \setRequiredRooms(c)} \sum_{p \in \setPeriods} \frac{\parSeatsInRoom(r) - \parSeatsForCourse(c)}{\parSeatsInRoom(r)} \cdot \varAssignment_{c,p,r} 
\end{align}

\noindent
Assignment of courses to periods that are generally unpopular with students, taking into account the sizes of the courses.
The factor $0.5 + (\parSeatsForCourse(c))/(2 \cdot \max_{c\in\setCourses}(\parSeatsForCourse(c)))$ is defined such that it is between 0.5 and 1 for each course~$c\in\setCourses$, depending on the expected number of students in~$c$. Larger courses are penalized more if they are assigned to such periods:
\begin{align}
z_{18}(x)=\sum_{p \in \setPeriods} \sum_{c \in \setCourses} \sum_{r \in \setRequiredRooms(c)} \parPeriodPenalty(p) \cdot \left( 0.5 + \frac{\parSeatsForCourse(c)}{2 \cdot \max_{c\in\setCourses}(\parSeatsForCourse(c))} \right) \cdot \varAssignment_{c,p,r} 
\end{align}

\noindent
Assignment of courses to rooms that should preferably remain unassigned to preserve them for independent student study:
\begin{align}
z_{19}(x)=\sum_{c \in \setCourses} \sum_{r \in \setRequiredRooms(c)} \sum_{p \in \setPeriods} \parRoomPenalty(r) \cdot \varAssignment_{c,p,r} 
\end{align}

\subsubsection{Perturbations Objective}\label{UTT:perturbation_objective}
The perturbations objective~$\parPerturbationsObj$ counts the number of perturbations, where each session assignment of a course~$c\in\setCourses$ with~$\setPeriodsOfLastYear(c)\neq\emptyset$ to a period outside~$\setPeriodsOfLastYear(c)$ counts as one perturbation:

\begin{align}
\parPerturbationsObj(x)=\sum_{c \in \setCourses : \setPeriodsOfLastYear(c)\neq\emptyset} \sum_{p \in \setPeriods\setminus\setPeriodsOfLastYear(c)} \sum_{r \in \setRequiredRooms(c)} \varAssignment_{c,p,r}
\end{align}

\subsection{Constraints}
In this subsection, we present the constraints of our MIP model. We distinguish between general constraints, lecturers’ constraints, and students’ constraints.

\subsubsection{General Constraints}\label{UTT:general_constraints}
\noindent
Each course \(c \in \setCourses\) must be scheduled in exactly~$\parTimesPerWeek(c)$ many sessions, each time in a room that satisfies the requirements of~$c$:
\begin{align}
\sum_{p \in \setPeriods} \sum_{r \in \setRequiredRooms(c)} \varAssignment_{c,p,r} &= \parTimesPerWeek(c), \quad \forall c \in \setCourses  
\end{align}

\noindent
No two rooms~$r_1$ and~$r_2$ with $\{r_1, r_2\} \in \setRoomConflicts$ can be used simultaneously. In particular, if a large room can be split into two smaller rooms, the corresponding subrooms cannot be used if the large room is active:
\begin{align}
&\varRoomAvailability_{r_1, p} + \varRoomAvailability_{r_2, p} \leq 1, \quad\forall\{r_1, r_2\}\in\setRoomConflicts, p \in \setPeriods
\end{align}

\noindent
For each room in $\setRoomDaySetting$, a single decision per day must be made whether the room is active throughout that day:
\begin{align}
&\varRoomAvailability_{r, p} - \varRoomAvailability_{r, p+1} = 0, \\ &\forall r \in \setRoomDaySetting, d \in \setDays, p \in \{\min(\setPeriodsAtDay(d)), \dots , \max(\setPeriodsAtDay(d))-1\}\notag
\end{align}

\noindent
A room in $\setRoomDaySetting$ must be active in every period in which a course is assigned to it:
 \begin{align}
&\sum_{c \in \setCompatibleCourses(r)} \varAssignment_{c,p,r} - \varRoomAvailability_{r, p} \leq 0,\quad\forall r \in \setRooms, p \in \setPeriods
\end{align}

\noindent
Rooms are not available during their excluded periods:
\begin{align}
&\sum_{r \in \setRooms}\sum_{p \in \setExcludedPeriodsOfRoom(r)}\sum_{c \in \setCompatibleCourses(r)} \varAssignment_{c,p,r} = 0
\end{align}

\subsubsection{Lecturers' Constraints}\label{UTT:lecturers_constraints}
\noindent
Lecturers cannot teach in their excluded periods:
\begin{align}
&\sum_{l \in \setLecturers}\sum_{p \in \setExcludedPeriods(l)}\sum_{c \in \setCoursesOfLecturer(l)} \sum_{r \in \setRequiredRooms(c)} \varAssignment_{c, p, r} = 0
\end{align}

\noindent
The variable~$\varLecAtPeriod_{l, p}$ is set to~1 if and only if lecturer~$l\in\setLecturers$ teaches in period~$p\in\setPeriods$:
\begin{align}
&\sum_{c \in \setCoursesOfLecturer(l)} \sum_{r \in \setRequiredRooms(c)} \varAssignment_{c,p,r} - \varLecAtPeriod_{l, p} = 0,\quad\forall l \in \setLecturers, p \in \setPeriods
\end{align}

\noindent
The variable~$\varLecAtDay_{l, d}$ is set to~1 if and only if lecturer~$l\in \setLecturers$ teaches on day~$d\in\setDays$:
\begin{align}
&\sum_{p \in \setPeriodsAtDay(d)} \sum_{c \in \setCoursesOfLecturer(l)} \sum_{r \in \setRequiredRooms(c)} \varAssignment_{c,p,r} - \lvert\setPeriodsAtDay(d)\rvert \cdot \varLecAtDay_{l, d} \leq 0,\quad
&\forall l \in \setLecturers, d \in \setDays\\
&\varLecAtDay_{l, d} - \sum_{p \in \setPeriodsAtDay(d)} \sum_{c \in \setCoursesOfLecturer(l)} \sum_{r \in \setRequiredRooms(c)} \varAssignment_{c,p,r} \leq 0,\quad
&\forall l \in L, d \in \setDays
\end{align}

\noindent
No two courses taught by the same lecturer can take place at the same time:
\begin{align}
\sum_{c \in \setCoursesOfLecturer(l)} \sum_{r \in \setRequiredRooms(c)} \varAssignment_{c,p,r} &\leq 1, \quad \forall l \in \setLecturers, p \in \setPeriods  
\end{align}

\noindent
Courses must take place in their predefined periods if this is specified by their lecturers:
\begin{align}
&\sum_{r \in \setRequiredRooms(c)} \varAssignment_{c,p,r} = 1,\quad
\forall c \in \setCourses, p \in \setfixPeriods(c)
\end{align}

\noindent
For each period group~$G \in \setPeriodGroups$, the number of periods in~$G$ in which a lecturer teaches must not exceed the specified limit~$\parMaxPeriodSet(G)$:
\begin{align}
&\sum_{p \in G} \sum_{c \in \setCoursesOfLecturer(l)} \sum_{r \in \setRequiredRooms(c)} \varAssignment_{c,p,r} \leq \parMaxPeriodSet(G),\quad
&\forall l \in \setLecturers, G \in \setPeriodGroups
\end{align}

\noindent
The variable~$\varLecConsecutivePeriods_{l,p}$ is set to~1 if lecturer~$l\in\setLecturersDistPeriods$ teaches in the two consecutive periods~$p\in\setSubsequentPeriods$ and~$p+1$ on the same day:
\begin{align}
&\varLecAtPeriod_{l, p} + \varLecAtPeriod_{l, p+1} - \varLecConsecutivePeriods_{l, p} \leq 1,\quad
& \forall l \in \setLecturersDistPeriods, p \in \setSubsequentPeriods
\end{align}

\noindent
The variable~$\varLecDisjointPeriods_{l,p}$ is set to~1 if lecturer~$l\in\setLecturersCompactPeriods$ teaches in period~$p$ and in a later period on the same day, with at least one free period in between:
\begin{align}
&\varLecAtPeriod_{l, p} - 
\varLecAtPeriod_{l, p + 1} + \varLecAtPeriod_{l, p + 2} - \varLecDisjointPeriods_{l, p} \leq 1, \\ 
&\forall l \in \setLecturersCompactPeriods, d \in \setDays, p \in \{\min(\setPeriodsAtDay(d)), \dots , \max(\setPeriodsAtDay(d)) -2 \} \notag\\
&\varLecAtPeriod_{l, p} - \varLecAtPeriod_{l, p+1} - \varLecAtPeriod_{l, p + 2} + \varLecAtPeriod_{l, p + 3} - \varLecDisjointPeriods_{l, p} \leq 1, \\ &\forall l \in \setLecturersCompactPeriods, d \in \setDays, p \in \{ \min(\setPeriodsAtDay(d)), \dots , \max(\setPeriodsAtDay(d)) -3 \} \notag \\
&\varLecAtPeriod_{l, p} - \varLecAtPeriod_{l, p + 1} - \varLecAtPeriod_{l, p + 2} - \varLecAtPeriod_{l, p + 3} + \varLecAtPeriod_{l, p + 4} - 
\varLecDisjointPeriods_{l, p} \leq 1, \\ &\forall l \in \setLecturersCompactPeriods, d \in \setDays, p \in \{ \min(\setPeriodsAtDay(d)), \dots , \max(\setPeriodsAtDay(d)) -4 \} \notag\\
&\varLecAtPeriod_{l, p} - \varLecAtPeriod_{l, p + 1} - 
\varLecAtPeriod_{l, p + 2} - 
\varLecAtPeriod_{l, p + 3} - 
\varLecAtPeriod_{l, p + 4} + 
\varLecAtPeriod_{l, p + 5} - \varLecDisjointPeriods_{l, p} \leq 1, \\ 
&\forall l \in \setLecturersCompactPeriods, d \in \setDays, p \in \{ \min(\setPeriodsAtDay(d)), \dots , \max(\setPeriodsAtDay(d)) -5 \} \notag
\end{align}

\noindent
The variable~$\varLecAtDayPlusOne_{l,d}$ is set to~1 if lecturer~$l\in\setLecturersDistWeek$ teaches on day~$d\in \{ 1, \dots, \lvert\setDays \rvert - 1 \}$ and on day~$d+1$:
\begin{align}
&\varLecAtDay_{l, d} + \varLecAtDay_{l, d+1} - \varLecAtDayPlusOne_{l, d} \leq 1,\quad\forall l \in \setLecturersDistWeek, d \in \{ 1, \dots, \lvert\setDays \rvert - 1 \}
\end{align}

\noindent
The variables~$\varLecAtDayPlusTwo_{l,d}$, $\varLecAtDayPlusThree_{l,d}$, and $\varLecAtDayPlusFour_{l,d}$ are set to~1 if lecturer~$l\in\setLecturersCompactWeek$ teaches on day~$d$ and again on day~$d+2$, $d+3$, or~$d+4$, respectively, with no teaching on the intervening day(s):
\begin{align}
&\varLecAtDay_{l, d} - \varLecAtDay_{l, d+1} + \varLecAtDay_{l, d+2} - \varLecAtDayPlusTwo_{l, d} \leq 1, \\
&\forall l \in \setLecturersCompactWeek, d \in \{ 1, \dots, \lvert \setDays \rvert - 2 \} \notag\\
&\varLecAtDay_{l, d} - \varLecAtDay_{l, d+1} - \varLecAtDay_{l, d+2} + \varLecAtDay_{l, d+3} - \varLecAtDayPlusThree_{l, d} \leq 1, \\
&\forall l \in \setLecturersCompactWeek, d \in \{ 1, \dots, \lvert \setDays \rvert - 3 \} \notag\\
&\varLecAtDay_{l, d} - \varLecAtDay_{l, d+1} - \varLecAtDay_{l, d+2} - \varLecAtDay_{l, d+3} + \varLecAtDay_{l, d+4} - \varLecAtDayPlusFour_{l, d} \leq 1, \\
&\forall l \in \setLecturersCompactWeek, d \in \{ 1, \dots, \lvert\setDays \rvert - 4 \} \notag
\end{align}

\noindent
The variable~$\varLecDevDaysPos_{l}$ represents the positive deviation from the desired number of teaching days of lecturer~$l\in\setLecturersPrefDaysPerWeek$, whereas~$\varLecDevDaysNeg_{l}$ captures the corresponding negative deviation:
\begin{align}
&\sum_{d \in \setDays} \varLecAtDay_{l, d} - \varLecDevDaysPos_{l} + \varLecDevDaysNeg_{l} = \parNumberOfDays(l),\quad\forall l \in \setLecturersPrefDaysPerWeek
\end{align}

\noindent
The total number of periods that a lecturer~$l\in\setLecturersPrefDaysPerWeek$ is required to teach is given by~$\sum_{c \in \setCoursesOfLecturer(l)} \parTimesPerWeek(c)$, while the desired number of teaching days is specified by~$\parNumberOfDays(l)$. To distribute the teaching load evenly across these days, each teaching day should ideally comprise~$\frac{\sum_{c \in \setCoursesOfLecturer(l)} \parTimesPerWeek(c)}{\parNumberOfDays(l)}$ teaching periods. For each teaching day~$d\in\setDays$, the variable~$\varLecDeviationCoursesPerDay_{l, d}$ captures the deviation from the corresponding rounded-up or rounded-down desired average number of teaching periods for lecturer~$l$:
\begin{align}
&\sum_{p \in \setPeriodsAtDay(d)} \sum_{c \in \setCoursesOfLecturer(l)} \sum_{r \in \setRequiredRooms(c)} \varAssignment_{c,p,r} - \left\lceil \frac{\sum_{c \in \setCoursesOfLecturer(l)} \parTimesPerWeek(c)}{\parNumberOfDays(l)} \right\rceil\leq\varLecDeviationCoursesPerDay_{l, d}
, \\
&\forall l \in \setLecturersPrefDaysPerWeek, d \in \setDays \notag\\
&\varLecAtDay_{l, d} \cdot \left\lfloor \frac{\sum_{c \in \setCoursesOfLecturer(l)} \parTimesPerWeek(c)}{\parNumberOfDays(l)} \right\rfloor - \sum_{p \in \setPeriodsAtDay(d)} \sum_{c \in \setCoursesOfLecturer(l)} \sum_{r \in \setRequiredRooms(c)} \varAssignment_{c,p,r} \leq \varLecDeviationCoursesPerDay_{l, d}, \\
&\forall l \in \setLecturersPrefDaysPerWeek, d \in \setDays \notag
\end{align}

\noindent
Courses of a lecturer that are not follow-up courses should take place in the same room. The variable~$\varConsCourseSameRoomVio_{l, p}$ is set to~1 if two such courses of lecturer~$l\in\setLecturers$ are scheduled consecutively in period~$p\in \setSubsequentPeriods$ and $p+1$, but in different rooms:
\begin{align}
&\sum_{c \in \setCoursesOfLecturer(l)} \, \sum_{r_1 \in \setRequiredRooms(c)} \varAssignment_{c, p, r_1} 
+ \sum_{c \in \setCoursesOfLecturer(l)} \sum_{r_2 \in \setRequiredRooms(c) \setminus \{r_1\}} \varAssignment_{c, p+1, r_2} - \varConsCourseSameRoomVio_{l, p} \leq 1,\\
&\forall l \in \setLecturers, p \in \setSubsequentPeriods \notag
\end{align}    

\subsubsection{Students' Constraints}\label{UTT:students_constraints}
\noindent
The variable~$\varCurriculumAtPeriod_{\kappa, p}$ is set to~1 if an elective lecture or exercise class of curriculum~$\kappa\in\setCurricula$ is assigned to period~$p\in\setPeriods$:
\begin{align}
&\sum_{r \in \setRequiredRooms(c)} \varAssignment_{c,p,r} - \varCurriculumAtPeriod_{\kappa, p} \leq 0,\quad\forall \kappa \in \setCurricula, c \in (\setLectures\cup\setExercises) \cap \setElectiveCourses(\kappa), p \in \setPeriods
\end{align}

\noindent
The~variable $\varCurriculumMandatoryAtPeriod_{\kappa, p}$ is set to~1 if a mandatory lecture or exercise class of curriculum~$\kappa\in\setCurricula$ is assigned to period~$p\in\setPeriods$:
\begin{align}
&\sum_{r \in \setRequiredRooms(c)} \varAssignment_{c,p,r} - \varCurriculumMandatoryAtPeriod_{\kappa, p} \leq 0,\quad\forall \kappa \in \setCurricula, c \in (\setLectures\cup\setExercises) \cap \setMandCourses(\kappa), p \in \setPeriods
\end{align}

\noindent
At most one mandatory lecture or exercise class in a curriculum can take place in a period:
\begin{align}
&\sum_{c \in \setMandCourses(\kappa)\cap(\setLectures\cup\setExercises)} \sum_{r \in \setRequiredRooms(c)} \varAssignment_{c, p, r} \leq 1,\quad\forall \kappa \in \setCurricula, p \in \setPeriods  
\end{align}

\noindent
At most one lecture or exercise class from the same module can take place in a period:
\begin{align}
&\sum_{c\in(\setLectures\cup\setExercises)\cap M} \sum_{r \in \setRequiredRooms(c)} \varAssignment_{c,p,r} \leq 1,\quad\forall M \in \setModules, p \in \setPeriods
\end{align}

\noindent
No tutorial can take place at the same time as a lecture or exercise class of the same module:
\begin{align}
&\sum_{r \in \setRequiredRooms(t)} \varAssignment_{t, p, r} + \parTimesPerWeek(t) \cdot \sum_{c \in (\setLectures \cup \setExercises ) \cap M} \sum_{r \in \setRequiredRooms(c)} \varAssignment_{c,p,r} \leq \parTimesPerWeek(t), \\
&\forall M \in \setModules, t \in \setTutorials \cap M, p \in \setPeriods \notag
\end{align}

\noindent
For each mandatory tutorial~$c\in\setMandCourses(\kappa)\cap\setTutorials$ of a curriculum~$\kappa\in\setCurricula$, exactly one period~$p\in\setPeriods$ with~$\varTutorium_{\kappa, c, p}=1$ must be selected such that students of curriculum~$\kappa$ can attend one session of~$c$:
\begin{align}
&\sum_{p \in \setPeriods} \varTutorium_{\kappa, c, p} = 1, \quad\forall \kappa \in \setCurricula, c \in \setMandCourses(\kappa) \cap \setTutorials\\
&\varTutorium_{\kappa, c, p} - \sum_{r \in \setRequiredRooms(c)} \varAssignment_{c,p,r} \leq 0, \quad\forall \kappa \in \setCurricula, c \in \setMandCourses(\kappa) \cap \setTutorials, p \in \setPeriods
\end{align}

\noindent
For each curriculum~$\kappa\in\setCurricula$ and each period~$p\in\setPeriods$, at most one mandatory tutorial session selected for students of curriculum~$\kappa$ can take place. Moreover, no mandatory lecture or exercise class of curriculum~$\kappa$ can take place in the same period as such a selected tutorial session:
\begin{align}
&\sum_{c \in \setMandCourses(\kappa) \cap \setTutorials}
\varTutorium_{\kappa, c, p} \leq 1,\quad\forall \kappa \in \setCurricula, \forall p \in \setPeriods\\
&\varTutorium_{\kappa, c_1, p} + \sum_{r \in \setRequiredRooms(c_2)} \varAssignment_{c_2,p,r}\leq 1, \\
&\forall \kappa \in \setCurricula, c_1 \in \setMandCourses(\kappa) \cap \setTutorials,\notag\\ &c_2 \in (\setLectures\cup\setExercises) \cap \setMandCourses(\kappa), p \in \setPeriods \notag
\end{align}

\noindent
The variable~$\varCourseDistVio_{c_1, c_2, d}$ is set to~1 if courses~$c_1, c_2\in\setCourses$ with~$\{c_1,c_2\}\in\setCourseTuples$ are both scheduled on day~$d\in\setDays$, even though they should preferably not be scheduled on the same day:
\begin{align}
&\sum_{p \in \setPeriodsAtDay(d)} \sum_{r \in \setRequiredRooms(c_1)} \varAssignment_{c_1,p,r}  + \sum_{p \in \setPeriodsAtDay(d)} \sum_{r \in \setRequiredRooms(c_2)}\varAssignment_{c_2,p,r} - \varCourseDistVio_{c_1, c_2, d} \leq 1, \\
&\forall \{c_1, c_2\} \in \setCourseTuples, d \in \setDays \notag
\end{align}

\noindent
Courses~$c_1,c_2\in\setCourses$ with $\{c_1, c_2\}\in\setCourseConflict$ must not be scheduled on the same day: 
\begin{align}
&\sum_{p_1 \in \setPeriodsAtDay(d)}\sum_{r_1 \in \setRequiredRooms(c_1)} \varAssignment_{c_1,p_1,r_1} + \sum_{p_2 \in \setPeriodsAtDay(d)}\sum_{r_2 \in \setRequiredRooms(c_2)} \varAssignment_{c_2,p_2,r_2}\leq 1,\quad\forall \{c_1,c_2\} \in \setCourseConflict, d\in\setDays
\end{align}

\noindent
The follow-up course~$\parFollowingCourse(c)$ of a course~$c\in\setCourses$ has to be scheduled immediately after~$c$ in the next time slot on the same day and in the same room: 
\begin{align}
&\sum_{p\in\setPeriods\setminus\setSubsequentPeriods}\sum_{r\in\setRequiredRooms(c)}\varAssignment_{c, p, r} = 0, \quad
\forall c \in \setCourses : \parFollowingCourse(c)\neq\texttt{None}\\
&\varAssignment_{c, p, r} - \varAssignment_{\parFollowingCourse(c), p+1, r} = 0, \\
&\forall c \in \setCourses : \parFollowingCourse(c)\neq\texttt{None}, r \in \setRequiredRooms(c) \cap \setRequiredRooms(\parFollowingCourse(c)), p\in \setSubsequentPeriods \notag\\
&\varAssignment_{c, p, r} = 0, \label{app:following_course_same_room_}\\
&\forall c \in \setCourses : \parFollowingCourse(c)\neq\texttt{None}, r \in \setRequiredRooms(c) \setminus \setRequiredRooms(\parFollowingCourse(c)), p \in \setPeriods \notag
\end{align}

\noindent
The courses in each module~$M\in\setModules$ should be scheduled according to the specified order over the week. If~$c_1\prec c_2$ for~$c_1,c_2\in M$, then the first session of~$c_1$ should be scheduled before all sessions of~$c_2$. If the sessions of those two courses are scheduled differently, the variable~$\varCourseOrderViolation_{c_1, c_2}$ is set to~1:
\begin{align}
&\sum_{p \in \setPeriods: p \leq i} \sum_{r \in \setRequiredRooms(c_2)} \varAssignment_{c_2, p, r} - \parTimesPerWeek(c_2)\cdot\sum_{p \in \setPeriods: p \leq i} \sum_{r \in \setRequiredRooms(c_1)} \varAssignment_{c_1, p, r} - \varCourseOrderViolation_{c_1, c_2} \leq 0,\\
&\forall M \in \setModules, c_1, c_2 \in M : c_1 \prec c_2, i \in \setPeriods \notag
\end{align}

\noindent
The variable~$\varCurCourseVio_{\kappa, c, p}$ is set to~1 if an elective lecture or exercise class~$c\in(\setLectures \cup \setExercises) \cap \setElectiveCourses(\kappa)$ of curriculum~$\kappa\in\setCurricula$ overlaps with another course of curriculum~$\kappa$ in period~$p\in\setPeriods$. The second constraint ensures that~$\varCurCourseVio_{\kappa, c, p}$ can only be set to~1 for a course~$c$ that is actually assigned to period~$p$:
\begin{align}
&\sum_{c \in (\setLectures \cup \setExercises ) \cap \setOfCoursesInCurr(\kappa) } \sum_{r \in \setRequiredRooms(c)} \varAssignment_{c, p, r} - \sum_{c \in (\setLectures \cup \setExercises) \cap \setElectiveCourses(\kappa) } \varCurCourseVio_{\kappa, c, p} \leq 1, \\
&\forall \kappa \in \setCurricula, p \in \setPeriods\notag\\
&\varCurCourseVio_{\kappa, c, p}  - \sum_{r \in \setRequiredRooms(c)} \varAssignment_{c, p, r}\leq 0,\quad\forall \kappa \in \setCurricula, c \in (\setLectures \cup \setExercises) \cap \setElectiveCourses(\kappa), p \in \setPeriods
\end{align}

\noindent
The variable~$\varCurMand_{\kappa, p}$ is set to~1 if two mandatory lectures or exercise classes of curriculum~$\kappa\in\setCurricula$ are assigned to periods~$p\in\setSubsequentPeriods$ and~$p+1$:
\begin{align}
&\varCurriculumMandatoryAtPeriod_{\kappa, p} + \varCurriculumMandatoryAtPeriod_{\kappa, p+1} - \varCurMand_{\kappa, p} \leq 1,\quad\forall \kappa \in \setCurricula, \forall p\in \setSubsequentPeriods
\end{align}

\noindent
The variable~$\varCurThreePeriods_{\kappa, p}$ is set to~1 if, in curriculum~$\kappa\in\setCurricula$, two mandatory and one elective lecture or exercise class are scheduled in period~$p\in\setPeriods$ and the two subsequent periods on the same day:
\begin{align}
&\varCurriculumMandatoryAtPeriod_{\kappa, p} + \varCurriculumMandatoryAtPeriod_{\kappa, p+1} + \varCurriculumAtPeriod_{\kappa, p+2} - \varCurThreePeriods_{\kappa, p} \leq 2, \\
&\forall \kappa \in \setCurricula, d \in \setDays, p \in \{ \min(\setPeriodsAtDay(d)), \dots , \max(\setPeriodsAtDay(d)) -2 \} \notag\\
&\varCurriculumMandatoryAtPeriod_{\kappa, p} + \varCurriculumAtPeriod_{\kappa, p+1} + \varCurriculumMandatoryAtPeriod_{\kappa, p+2} - \varCurThreePeriods_{\kappa, p} \leq 2, \\
&\forall \kappa \in \setCurricula, d \in \setDays, p \in \{ \min(\setPeriodsAtDay(d)), \dots , \max(\setPeriodsAtDay(d)) -2 \} 
\notag\\
&\varCurriculumAtPeriod_{\kappa, p} + \varCurriculumMandatoryAtPeriod_{\kappa, p+1} + \varCurriculumMandatoryAtPeriod_{\kappa, p+2} - \varCurThreePeriods_{\kappa, p} \leq 2, \\
&\forall \kappa \in \setCurricula, d \in \setDays, p \in \{ \min(\setPeriodsAtDay(d)), \dots , \max(\setPeriodsAtDay(d)) -2 \}
\notag 
\end{align}

\noindent
The variable~$\varCurThreeMandPeriods_{\kappa, p}$ is set to~1 if, in curriculum~$\kappa\in\setCurricula$, three mandatory lectures or exercise classes are scheduled in period~$p\in\setPeriods$ and the two subsequent periods on the same day:
\begin{align}
&\varCurriculumMandatoryAtPeriod_{\kappa, p} + \varCurriculumMandatoryAtPeriod_{\kappa, p+1} + \varCurriculumMandatoryAtPeriod_{\kappa, p+2} - \varCurThreeMandPeriods_{\kappa, p} \leq 2, \label{UTT:last_constraint}\\
&\forall \kappa \in \setCurricula, d \in \setDays, p \in \{ \min(\setPeriodsAtDay(d)), \dots , \max(\setPeriodsAtDay(d)) -2 \} \notag
\end{align} \section{Multi-Objective Optimization Approach}\label{UTT:sec:methodology}

In this section, we introduce key concepts of multi-objective optimization tailored to our problem setting (Section~\ref{UTT:subsec:multiobj}) and subsequently describe the multi-objective solution approach we propose (Section~\ref{UTT:subsec:algorithm}).

\subsection{Multi-Objective Optimization Concepts}\label{UTT:subsec:multiobj}

To analyze the trade-offs between the students' and lecturers' objectives, we briefly introduce the concepts from multi-objective optimization required in the following. The definitions below are tailored to the bi-objective setting considered in this paper. For more general definitions, see~\citet{ehrgott_multicriteria_2005}. Recall that~$x$ denotes the vector comprising all variables of the MIP model introduced in Section~\ref{UTT:sec:MIP}.

\begin{definition}
Let~$\setSolutions$ denote the set of all feasible solutions~$x$ of the MIP model. For~$x \in \setSolutions$, the vector
$z(x) \coloneqq \big(\parLecturersObj(x), \parStudentsObj(x)\big)\in\mathbb{R}^2$ is called  the \emph{image} of~$x$. The set $\setImage = \{ z(x) \mid x \in \setSolutions \}$ of all images of feasible solutions is called the \emph{image set},
which is a subset of the \emph{objective space}~$\mathbb{R}^2$.

A solution~$x \in \setSolutions$ \emph{dominates} a solution~$\tilde{x} \in \setSolutions$ if
$\parLecturersObj(x) \leq \parLecturersObj(\tilde{x})$,
$\parStudentsObj(x) \leq \parStudentsObj(\tilde{x})$,
and at least one of the two inequalities is strict. A solution~$x \in \setSolutions$ is \emph{efficient} if it is not dominated by any other solution in~$\setSolutions$. In this case, we call the corresponding image~$z(x)\in \setImage$ \emph{nondominated}.
The set of all nondominated images is denoted by
$\setNondomImages = \{ z(x) \mid x \in \setSolutions \text{ is efficient} \}$.
\end{definition}

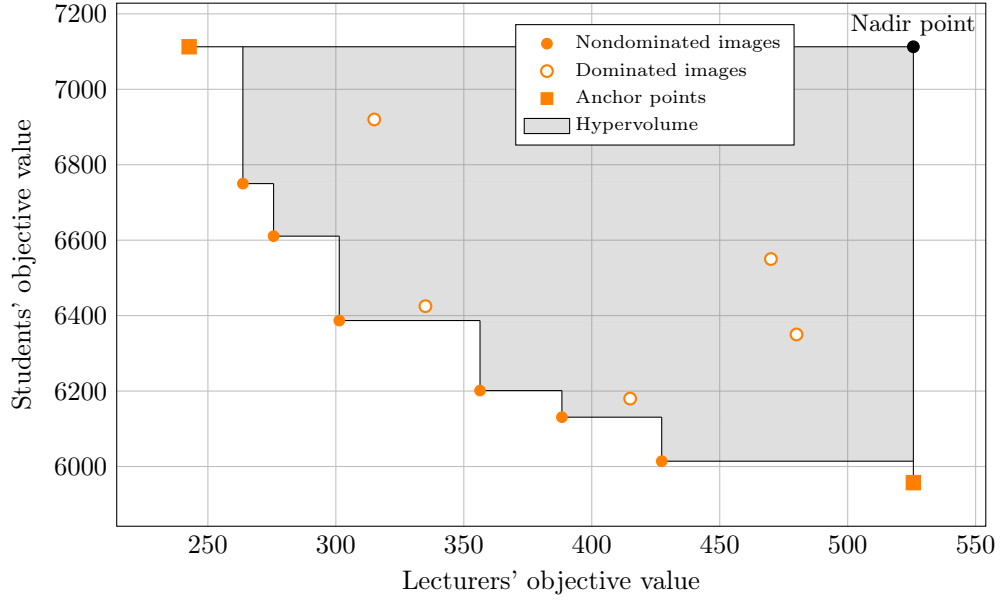
\begin{figure}[htbp]
\begin{tikzpicture}

\pgfplotsset{
    hvline/.style={color=black},
    hvfill/.style={fill=gray, fill opacity=0.25, draw=none},
    nadir/.style={only marks, mark=*, mark size=2.2pt, color=black}, 
    anchorpt/.style={only marks, mark=square*, mark size=2.5pt, thick},
    dominated/.style={only marks, mark=*, mark size=2.2pt,
                      mark options={fill=white, draw=orange, line width=0.8pt}}}
    
\begin{axis}[
    grid=both,
    xlabel={Lecturers' objective value},
    ylabel={Students' objective value},
    width=\linewidth,
    height=0.65\linewidth,
    legend style={at={(0.78,0.96)},anchor=north east, font=\footnotesize},
    legend cell align={left},
    scaled ticks=false,
    tick label style={/pgf/number format/1000 sep={}},
    tick style={draw=none},
    ]

\addplot[nadir, forget plot] coordinates {(525.6666667,7112.719918)};
    \node[font=\small] at (525.6666667,7170) {Nadir point};
    
\addplot[
      only marks,
      mark=*,
      color=orange
    ] coordinates {
          (525.6666667,5957.19)
          (427.3333333,6014.03)
          (388.3333333,6130.75)
          (356.3333333,6201.22)
          (301.3333333,6386.9)
          (275.6666667,6610.679995)
          (263.6666667,6749.98)
          (242.6666667,7112.719918)
    };
    \addlegendentry{Nondominated images}

    \addlegendimage{dominated}
    \addlegendentry{Dominated images}

    \addplot[anchorpt, color=orange, forget plot] coordinates {
    (525.6666667,5957.19)
    (242.6666667,7112.719918)};
    
    \addlegendimage{only marks, mark=square*, mark size=2.0pt, thick, color=orange}
    \addlegendentry{Anchor points}        
    \addlegendimage{area legend, fill=gray, fill opacity=0.25, draw=none}
    \addlegendentry{Hypervolume}

    \addplot[hvfill] coordinates {
       (525.6666667,6749.98)
       (263.6666667,6749.98)
       (263.6666667,7112.719918)
       (525.6666667,7112.719918)
    } \closedcycle;
    
    \addplot[hvfill] coordinates {
       (525.6666667,6610.679995)
       (275.6666667,6610.679995)
       (275.6666667,6749.98)
       (525.6666667,6749.98)
    } \closedcycle;
    
    \addplot[hvfill] coordinates {
       (525.6666667,6386.9)
       (301.3333333,6386.9)
       (301.3333333,6610.679995)
       (525.6666667,6610.679995)
    } \closedcycle;
    
    \addplot[hvfill] coordinates {
       (525.6666667,6201.22)
       (356.3333333,6201.22)
       (356.3333333,6386.9)
       (525.6666667,6386.9)
    } \closedcycle;
    
    \addplot[hvfill] coordinates {
       (525.6666667,6130.75)
       (388.3333333,6130.75)
       (388.3333333,6201.22)
       (525.6666667,6201.22)
    } \closedcycle;
    
    \addplot[hvfill] coordinates {
       (525.6666667,6014.03)
       (427.3333333,6014.03)
       (427.3333333,6130.75)
       (525.6666667,6130.75)
    } \closedcycle;

    \addplot[hvline] coordinates {(242.6666667,7112.719918) (525.6666667,7112.719918)};
    \addplot[hvline] coordinates {(525.6666667,5957.19) (525.6666667,7112.719918)};
    \addplot[hvline] coordinates {(263.6666667,6749.98) (263.6666667,7112.719918)};
    \addplot[hvline] coordinates {(263.6666667,6749.98) (275.6666667,6749.98)};
    \addplot[hvline] coordinates {(275.6666667,6610.679995) (275.6666667,6749.98)};
    \addplot[hvline] coordinates {(275.6666667,6610.679995) (301.3333333,6610.679995)};
    \addplot[hvline] coordinates {(301.3333333,6386.9) (301.3333333,6610.679995)};
    \addplot[hvline] coordinates {(301.3333333,6386.9) (356.3333333,6386.9)};
    \addplot[hvline] coordinates {(356.3333333,6201.22) (356.3333333,6386.9)};
    \addplot[hvline] coordinates {(356.3333333,6201.22) (388.3333333,6201.22)};
    \addplot[hvline] coordinates {(388.3333333,6130.75) (388.3333333,6201.22)};
    \addplot[hvline] coordinates {(388.3333333,6130.75) (427.3333333,6130.75)};
    \addplot[hvline] coordinates {(427.3333333,6014.03) (427.3333333,6130.75)};
    \addplot[hvline] coordinates {(427.3333333,6014.03) (525.6666667,6014.03)};

\addplot[dominated, forget plot] coordinates {
       (335,6425)
       (415,6180)
       (315,6920)
       (470,6550)
       (480,6350)
    };
        
\end{axis}
\end{tikzpicture}
\caption{The filled orange markers show nondominated images, i.e., images not dominated by any other image in the depicted set; among them, the square markers denote the anchor points. The unfilled orange markers show dominated images. The Nadir point is depicted in black. The area of the shaded gray region is the hypervolume dominated by the depicted images with respect to the Nadir point.
}\label{UTT:fig:hypervolume_example}
\end{figure}

In bi-objective optimization, the \emph{anchor points}, also called \emph{lexicographic points}, are the nondominated images obtained by optimizing one objective and then the other under the constraint that the first objective attains its optimal value. They correspond to the upper-left and lower-right corner points of the rectangle in~$\mathbb{R}^2$ that contains all nondominated images 
(see Figure~\ref{UTT:fig:hypervolume_example} for an illustration). In our setting, these points are formally defined as follows:

\begin{definition}
Let
$z^{\mathrm{lec}}_{\min} \coloneqq \min_{x \in \setSolutions} \parLecturersObj(x)$.
Then the \emph{lecturers' anchor point} is defined as
$\parLecAnchorPt \coloneqq z(x)$ for some
$x \in \arg\min \{ \parStudentsObj(x) \mid x \in \setSolutions,\; \parLecturersObj(x) = z^{\mathrm{lec}}_{\min} \}$.
Analogously, let
$z^{\mathrm{stud}}_{\min} \coloneqq \min_{x \in \setSolutions} \parStudentsObj(x)$.
Then the \emph{students' anchor point} is defined as
$\parStudAnchorPt \coloneqq z(x)$ for some
$x \in \arg\min \{ \parLecturersObj(x) \mid x \in \setSolutions,\; \parStudentsObj(x) = z^{\mathrm{stud}}_{\min} \}$.
\end{definition}

Computing the set~$\setNondomImages$ of nondominated images 
exactly is, in general, computationally challenging: The cardinality of~$\setNondomImages$ can be exponential in the input size, and deciding whether a given image is nondominated is NP-hard for many multi-objective optimization problems. 
The anchor points and the images computed by our algorithm are therefore approximate rather than exact, as the underlying optimization runs are terminated based on time limits or prescribed MIP optimality gaps. Accordingly, whenever we refer to anchor points, nondominated images, or corresponding solutions in the following, we mean the respective approximate images or solutions obtained in this way. 

The following definition introduces the \emph{(dominated) hypervolume} of a set of images, which is the most widely used quality indicator for sets of images in multi-objective optimization \citep[see, e.g.,][]{li2019quality} and will be used to compare different image sets obtained for the same instance in our computational experiments (see Figure~\ref{UTT:fig:hypervolume_example} for an illustration). A larger hypervolume indicates a better image set, as it corresponds to a larger dominated region in the objective space. The hypervolume is computed with respect to a reference point, which is commonly chosen as the \emph{Nadir point}, given by the componentwise worst objective values among all nondominated images, or as a slight modification of the Nadir point (e.g., $1.1$~times the Nadir point in the case of a minimization problem with nonnegative objectives). In our computational experiments, we use, for each instance, the Nadir point associated with the union of all image sets obtained under the considered bounds on the number of allowed perturbations. Accordingly, the following definition introduces the Nadir point for an arbitrary set of images in the objective space and defines the hypervolume with respect to an arbitrary reference point.

\begin{definition}\label{UTT:def:hypervolume}
Let $\setApproxSet \subseteq \setImage$ be a set of images. The \emph{Nadir point associated with~$\setApproxSet$} is defined as
$\parApproxNadirPt \coloneqq \left( \max_{z \in \setApproxSet} z_1,\; \max_{z \in \setApproxSet} z_2 \right)$.
If~$\parRefPt\in\mathbb{R}^2$ is a reference point satisfying $z_i\leq(\parRefPt)_i$ for all $z\in\setApproxSet$ and $i\in\{1,2\}$, the \emph{hypervolume} induced by~$\setApproxSet$ with respect to~$\parRefPt$ is defined as the Lebesgue measure of the set
\begin{align*}
\bigcup_{z \in \setApproxSet} \left[ z_1, (\parRefPt)_1 \right] \times \left[ z_2, (\parRefPt)_2 \right].
\end{align*}
\end{definition}

Note that the hypervolume of a set~$\setApproxSet \subseteq \setImage$ of images with respect to a given reference point as in Definition~\ref{UTT:def:hypervolume} can be computed easily in~$\mathcal{O}(|\setApproxSet|)$ time~\citep{kuhn2016hypervolume}.

\subsection{Multi-Objective Optimization Algorithm}\label{UTT:subsec:algorithm}
In this section, we describe our algorithm based on the MIP model from Section~\ref{UTT:sec:MIP}. The pseudocode is provided in Algorithm~\ref{UTT:algorithm}, and the pseudocode of the subroutines used within the algorithm is given separately in Subroutines~\ref{UTT:subroutine_lec_anchor_point}--\ref{UTT:subroutine_stud_e-constraint}.

We first determine the minimum required number~$\omega$ of perturbations for which a feasible timetable exists by minimizing~$\parPerturbationsObj$. This value is used to define the set~$\Theta = \{\omega, \omega+2, \omega+4, 20, 25, 30\}$ of allowed numbers of perturbations, chosen to cover different regions of the objective space. These values are tailored to our instances, for which no more than~$30$ perturbations were considered acceptable in practice, and may be adapted for other instances. For each value~$\theta\in\Theta$, we then use the lexicographic $\varepsilon$-constraint method~\citep[see, e.g.,][]{hamacher_finding_2007} with respect to the lecturers' and students' objectives to compute a set of images covering different parts of the objective space and a corresponding set of solutions, while bounding~$\parPerturbationsObj$ from above by~$\theta$. To this end, we solve a sequence of MIP models with varying upper bounds on either the lecturers' or the students' objective. Each solution is computed lexicographically by first optimizing one of these two objectives and then the other. However, since each optimization run is terminated after a predefined time limit of eight hours or once a MIP gap of 5\% is reached, optimality cannot always be guaranteed. Consequently, not all of the computed images are guaranteed to be nondominated.

In each optimization run, one of the three objectives, $\parLecturersObj(x)$, $\parStudentsObj(x)$, or $\parPerturbationsObj(x)$, is minimized, while the remaining objectives may be bounded from above. More precisely, we impose the bounds
\begin{align*}
\parPerturbationsObj(x) \leq \parPerturbationsBound,\qquad
\parLecturersObj(x) \leq \parLecturersBound,\qquad
\parStudentsObj(x) \leq \parStudentsBound,\qquad
\end{align*}
where each of the upper bounds $\parPerturbationsBound$, $\parLecturersBound$, and $\parStudentsBound$ may also take the value $+\infty$, and we assign different values to them in our algorithm. The full set of constraints from Section~\ref{UTT:sec:MIP} is included in every optimization run. 
In the following description, the notation
\begin{align*}
    (\hat{z}, \parPerturbationsBound, \parLecturersBound, \parStudentsBound)\rightarrow x
\end{align*}
indicates that the model is optimized with respect to the objective function~$\hat{z}\in\{\parLecturersObj,\parStudentsObj,\parPerturbationsObj\}$, using the specified values of the perturbations, lecturers', and students' bounds. The symbol~\(x\) on the right-hand side denotes the solution obtained from the corresponding optimization run.

\smallskip

We now describe the individual steps of the algorithm. The line numbers in square brackets refer to Algorithm~\ref{UTT:algorithm}. Throughout the algorithm, we maintain a set~$\setApproxSolutions$ of computed solutions, initialized as the empty set in line~\ref{UTT:alg:Solution_set}.

\begin{algorithm}[tbp]
\caption{Computation of Representative Solutions Under Perturbation Constraints}
\begin{algorithmic}[1]
\State Initialize $\setApproxSolutions\coloneqq\emptyset$ \label{UTT:alg:Solution_set}
\State Determine the minimum required number~$\omega$ of perturbations as $\omega \coloneqq \parPerturbationsObj(x)$ for $(\parPerturbationsObj, +\infty, +\infty, +\infty)\rightarrow x$\label{UTT:alg:det_min_feasible_pert}
\State Define set $\Theta = \{\omega, \omega+2, \omega+4, 20, 25, 30\}$ of allowed perturbations\label{UTT:alg:Set_allowed_pert}
\For{each number~$\theta \in \Theta$ of allowed perturbations}\label{UTT:alg:for_allowed_pert}
    \State Compute the lecturers' anchor point~$\parLecAnchorPtTheta$ using Subroutine~\ref{UTT:subroutine_lec_anchor_point}
    \State Compute the students' anchor point $\parStudAnchorPtTheta$ using Subroutine~\ref{UTT:subroutine_stud_anchor_point}
\EndFor\label{UTT:alg:for_allowed_pert_end}

\For{each number~$\theta \in \Theta$ of allowed perturbations}\label{UTT:alg:for_allowed_pert_2}
    \State Compute the lecturers' $\varepsilon$-constraint solutions using Subroutine~\ref{UTT:subroutine_lec_e-constraint}
    \State Compute the students' $\varepsilon$-constraint solutions using Subroutine~\ref{UTT:subroutine_stud_e-constraint}
\EndFor\label{UTT:alg:for_allowed_pert_2_end}

\State Return $\setApproxSolutions$ and $\setApproxSet = z(\setApproxSolutions)$\label{UTT:alg:return}
\end{algorithmic}\label{UTT:algorithm}
\end{algorithm}

\makeatletter
\renewcommand{\ALG@name}{Subroutine}
\makeatother
\setcounter{algorithm}{0}
\begin{algorithm}[tbp]
\caption{Compute lecturers' anchor point $\parLecAnchorPtTheta$}
\begin{algorithmic}[1]
\State Minimize $(\parLecturersObj, \theta, +\infty, +\infty)\rightarrow x$
\State Set $\parLecturersBound \coloneqq \parLecturersObj(x)$
\State Minimize $(\parStudentsObj, \theta, \parLecturersBound, +\infty)\rightarrow x$\label{UTT:alg:lec_anchor_opt_2}
\State Store the image of the solution as the lecturers' anchor point~$\parLecAnchorPtTheta = z(x)$\label{UTT:alg:store_lec_anchor}
\State Add~$x$ to~$\setApproxSolutions$\label{UTT:alg:store_lec_anchor2}
\end{algorithmic}\label{UTT:subroutine_lec_anchor_point}
\end{algorithm}

\begin{algorithm}[tbp]
\caption{Compute students' anchor point $\parStudAnchorPtTheta$}
\begin{algorithmic}[1]
\State Minimize $(\parStudentsObj, \theta, +\infty, +\infty)\rightarrow x$
\State Set $\parStudentsBound \coloneqq \parStudentsObj(x)$
\State Minimize $(\parLecturersObj, \theta, +\infty, \parStudentsBound)\rightarrow x$\label{UTT:alg:stud_anchor_opt_2}
\State Store the image of the solution as the students' anchor point~$\parStudAnchorPtTheta = z(x)$\label{UTT:alg:store_stud_anchor}
\State Add~$x$ to~$\setApproxSolutions$\label{UTT:alg:store_lec_anchor3}
\end{algorithmic}\label{UTT:subroutine_stud_anchor_point}
\end{algorithm}

\begin{algorithm}[tbp]
\caption{Compute lecturers' $\varepsilon$-constraint solutions}
\begin{algorithmic}[1]
\State $I_L \coloneqq [(\parLecAnchorPtTheta)_1,(\parStudAnchorPtTheta)_1]$
\State $\delta_L \coloneqq (\max(I_L)-\min(I_L))/4$
\State $\mathcal{E}_L \coloneqq \{\min(I_L)+\delta_L,\ \min(I_L)+\tfrac{5}{3}\delta_L,\ \min(I_L)+\tfrac{7}{3}\delta_L,\ \min(I_L)+3\delta_L\}$
\For{each $\varepsilon \in \mathcal{E}_L$}\label{UTT:alg:for_each_EL}
    \State Set $\parLecturersBound \coloneqq \varepsilon$\label{UTT:alg:set_eps_1}
    \State Minimize $(\parStudentsObj, \theta, \parLecturersBound, +\infty)\rightarrow x$\label{UTT:alg:min_lec_1}
    \State Set $\parStudentsBound \coloneqq \parStudentsObj(x)$\label{UTT:alg:set_stud_bound}
    \State Minimize $(\parLecturersObj, \theta, \parLecturersBound, \parStudentsBound)\rightarrow x$\label{UTT:alg:min_lec_2}
    \State Add~$x$ to~$\setApproxSolutions$\label{UTT:alg:store_lec_sol}
\EndFor
\end{algorithmic}\label{UTT:subroutine_lec_e-constraint}
\end{algorithm}

\begin{algorithm}[tbp]
\caption{Compute students' $\varepsilon$-constraint solutions}
\begin{algorithmic}[1]
\State $I_S \coloneqq [(\parStudAnchorPtTheta)_2,(\parLecAnchorPtTheta)_2]$
\State $\delta_S \coloneqq (\max(I_S)-\min(I_S))/4$
\State $\mathcal{E}_S \coloneqq \{\min(I_S)+\delta_S,\ \min(I_S)+\tfrac{5}{3}\delta_S,\ \min(I_S)+\tfrac{7}{3}\delta_S,\ \min(I_S)+3\delta_S\}$
\For{each $\varepsilon \in \mathcal{E}_S$}\label{UTT:alg:for_each_ES}
    \State Set $\parStudentsBound \coloneqq \varepsilon$\label{UTT:alg:set_eps_2}
    \State Minimize $(\parLecturersObj, \theta, +\infty, \parStudentsBound)\rightarrow x$\label{UTT:alg:min_stud_1}
    \State Set $\parLecturersBound \coloneqq \parLecturersObj(x)$\label{UTT:alg:set_lec_bound}
    \State Minimize $(\parStudentsObj, \theta, \parLecturersBound, \parStudentsBound)\rightarrow x$\label{UTT:alg:min_stud_2}
    \State Add~$x$ to~$\setApproxSolutions$\label{UTT:alg:store_stud_sol}
\EndFor
\end{algorithmic}\label{UTT:subroutine_stud_e-constraint}
\end{algorithm}
\makeatletter
\renewcommand{\ALG@name}{Algorithm}
\makeatother

\paragraph{Step 1: Determine the minimum number of perturbations [line~\ref{UTT:alg:det_min_feasible_pert}]}
We first compute the minimum required number~$\omega$ of perturbations by solving
\begin{align*}
    (\parPerturbationsObj, +\infty, +\infty, +\infty)\rightarrow x
\end{align*}
and setting~$\omega=\parPerturbationsObj(x)$. Consequently, no feasible solution with fewer than~$\omega$ perturbations exists. 

\paragraph{Step 2: Define the set of allowed perturbations [line~\ref{UTT:alg:Set_allowed_pert}]}
Using~$\omega$, we define the set of allowed numbers of perturbations as
\begin{align*}
\Theta = \{\omega,\omega + 2,\omega + 4,20,25,30\}.
\end{align*}
This choice reflects that the minimum number of perturbations observed in our instances ranges between eight and sixteen, while the planners consider at most 30 perturbations acceptable. For other instances, the set~$\Theta$ can be adapted accordingly.

\paragraph{Step 3: Compute anchor points [lines~\ref{UTT:alg:for_allowed_pert}--\ref{UTT:alg:for_allowed_pert_end}]}
For each value~$\theta \in \Theta$ for the allowed number of perturbations, we set~$\parPerturbationsBound=\theta$. We then compute the two anchor points defining the relevant range of lecturers' and students' objective values for the subsequent runs of the lexicographic $\varepsilon$-constraint method. We denote the lecturers' and students' anchor points obtained for this value of~$\theta$ by~$\parLecAnchorPtTheta$ and~$\parStudAnchorPtTheta$, respectively.

\begin{itemize}
    \item \emph{\textbf{Compute lecturers' anchor point $\parLecAnchorPtTheta$ (Subroutine~\ref{UTT:subroutine_lec_anchor_point})}}\\
    We first minimize $(\parLecturersObj, \theta, +\infty, +\infty)\rightarrow x$. We then set~$\parLecturersBound = \parLecturersObj(x)$ and perform a second optimization run with respect to the students' objective, $(\parStudentsObj, \theta, \parLecturersBound, +\infty)\rightarrow x$.
    The resulting solution~$x$ is added to~$\setApproxSolutions$, and its image~$z(x)\in \mathbb{R}^2$ defines the lecturers' anchor point~$\parLecAnchorPtTheta \coloneqq z(x)$.
    \smallskip
    \item \emph{\textbf{Compute students' anchor point $\parStudAnchorPtTheta$ (Subroutine~\ref{UTT:subroutine_stud_anchor_point})}}\\
    We first minimize $(\parStudentsObj, \theta, +\infty, +\infty)\rightarrow x$. We then set~$\parStudentsBound = \parStudentsObj(x)$ and perform a second optimization run with respect to the lecturers' objective, $(\parLecturersObj, \theta, +\infty, \parStudentsBound)\rightarrow x$.
    The resulting solution~$x$ is added to~$\setApproxSolutions$, and its image~$z(x)\in \mathbb{R}^2$ defines the students' anchor point~$\parStudAnchorPtTheta \coloneqq z(x)$.
\end{itemize}

\paragraph{Step 4: Generate $\varepsilon$-constrained solutions [lines~\ref{UTT:alg:for_allowed_pert_2}--\ref{UTT:alg:for_allowed_pert_2_end}]}
For each value~$\theta \in \Theta$ for the allowed number of perturbations, we use the lexicographic $\varepsilon$-constraint method to generate additional images between the two anchor points~$\parLecAnchorPtTheta$ and~$\parStudAnchorPtTheta$, together with corresponding solutions. Based on the lecturers' and students' objective values of the anchor points, we define four $\varepsilon$-values for each of the two objectives, which are used as upper bounds on the respective objective within the lexicographic $\varepsilon$-constraint method.

\begin{itemize}
    \item \emph{\textbf{Compute lecturers' $\varepsilon$-constraint solutions (Subroutine~\ref{UTT:subroutine_lec_e-constraint})}}\\
    We consider the interval $I_L \coloneqq [(\parLecAnchorPtTheta)_1,(\parStudAnchorPtTheta)_1]$ defined by the lecturers' objective values of the two anchor points and set $\delta_L \coloneqq (\max(I_L)-\min(I_L))/4$. Based on this interval, we define four $\varepsilon$-values, chosen more densely near the center of the interval, where the most balanced trade-offs between the two objectives are typically located:
    \begin{align*}
        \mathcal{E}_L \coloneqq
        \Big\{
        \min(I_L) + \delta_L,\;
        \min(I_L) + \tfrac{5}{3}\delta_L,\;
        \min(I_L) + \tfrac{7}{3}\delta_L,\;
        \min(I_L) + 3\delta_L
        \Big\}.
    \end{align*}
    These values are used as upper bounds on the lecturers' objective. For each~$\varepsilon \in \mathcal{E}_L$, we set~$\parLecturersBound\coloneqq\varepsilon$ and first minimize the students' objective, $(\parStudentsObj, \theta, \parLecturersBound, +\infty)\rightarrow x$. For the solution~$x$ obtained in this run, we then set~$\parStudentsBound\coloneqq\parStudentsObj(x)$ and perform a second optimization run with respect to the lecturers' objective, $(\parLecturersObj, \theta, \parLecturersBound, \parStudentsBound)\rightarrow x$. The resulting solution is added to~$\setApproxSolutions$ before the next value of~$\varepsilon$ is considered.
    \smallskip

        \item \emph{\textbf{Compute students' $\varepsilon$-constraint solutions (Subroutine~\ref{UTT:subroutine_stud_e-constraint})}}\\
    Analogously, we consider the interval $I_S \coloneqq [(\parStudAnchorPtTheta)_2,(\parLecAnchorPtTheta)_2]$ defined by the students' objective values of the two anchor points and set $\delta_S \coloneqq (\max(I_S)-\min(I_S))/4$. Based on this interval, we define four $\varepsilon$-values, again chosen more densely near the center of the interval, where the most balanced trade-offs between the two objectives are typically located:
    \begin{align*}
        \mathcal{E}_S \coloneqq
        \Big\{
        \min(I_S) + \delta_S,\;
        \min(I_S) + \tfrac{5}{3}\delta_S,\;
        \min(I_S) + \tfrac{7}{3}\delta_S,\;
        \min(I_S) + 3\delta_S
        \Big\}.
    \end{align*}
    These values are used as upper bounds on the students' objective. For each~$\varepsilon \in \mathcal{E}_S$, we set~$\parStudentsBound\coloneqq\varepsilon$ and first minimize the lecturers' objective, $(\parLecturersObj, \theta, +\infty, \parStudentsBound)\rightarrow x$. For the solution~$x$ obtained in this run, we then set~$\parLecturersBound\coloneqq\parLecturersObj(x)$ and perform a second optimization run with respect to the students' objective, $(\parStudentsObj, \theta, \parLecturersBound, \parStudentsBound)\rightarrow x$. The resulting solution is added to~$\setApproxSolutions$ before the next value of~$\varepsilon$ is considered.
\end{itemize}
The set of images computed by the algorithm is then given by
$\setApproxSet \coloneqq z(\setApproxSolutions)$. To improve the performance of subsequent optimization runs, we use the computed solutions as warm starts: whenever a new optimization run is initiated, the best solution in~$\setApproxSolutions$ satisfying all currently active bounds is provided as an initial solution. Hence, for each instance and each value~$\theta \in \Theta$ for the allowed number of perturbations, the algorithm computes up to ten solutions: two anchor-point solutions and eight additional $\varepsilon$-constraint solutions. Their images cover the relevant range of lecturers' and students' objective values between the two anchor points. \newpage
\section{Computational Experiments}\label{UTT:sec:computational_results}

In this section, we present our computational results obtained using the MIP formulation introduced in Section~\ref{UTT:sec:MIP} and the algorithm described in Section~\ref{UTT:subsec:algorithm}. The experiments are conducted on three real-world data sets from the TUM Campus Straubing from the winter semester 2024/2025~(Instance~1), the summer semester~2025~(Instance~2), and the winter semester 2025/2026~(Instance~3).

For each instance, we compute solution sets that represent different trade-offs between the students' and lecturers' objectives, using the algorithm described in Section~\ref{UTT:subsec:algorithm}. The aim is to analyze whether the objectives of the two stakeholder groups are conflicting and how strongly the bound on the number of allowed perturbations affects the resulting timetable quality. After consultation with the responsible planners, we set the maximum number of allowed perturbations to 30, which reflects the level of deviation from the previous time period assignments that they consider acceptable. Only for Instance~3, we additionally consider an unbounded number of allowed perturbations, since the corresponding running times remain moderate, as shown in Appendix~\ref{UTT:subsec:detailed_opt_results}. To compare the quality of the computed solution sets across different perturbation bounds, we evaluate their hypervolume in the objective space using a common reference point, namely the Nadir point associated with the union of all image sets computed for the considered perturbation bounds.

The computations were performed on a compute server with Ubuntu 20.04.6 LTS (GNU/Linux 5.4.0-208-generic x86\_64), an AMD EPYC 7542 processor (32 cores at 2.9~GHz, 64 threads), and 500~GB memory. 
All MIP models were solved with Gurobi 12.0.3 using a MIP gap of~5\% and a time limit of eight hours per optimization run.

\subsection{Instances}\label{UTT:subsec:instances}

The computational experiments are based on three real-world instances from the TUM Campus Straubing, comprising 162--181 courses, 80--95 lecturers, and 27~rooms; see Table~\ref{UTT:table:instances}. Most courses require only one session per week. The complete dataset for each semester was compiled by the academic planners, who provided the course information and detailed data on the available rooms. Lecturers specified their excluded periods, preferred periods, preferred number of teaching days, preferences regarding consecutive or non-consecutive teaching periods within a day, and preferences regarding consecutive or separated teaching days. They also provided course-specific information, including requirements for room assignments.

For the winter semester 2024/2025 (Instance~1), no information was available on the courses for which continuity of teaching periods, i.e., assignment to the same periods as in the corresponding previous winter semester, was desired. For this instance, we therefore randomly selected half of the courses as courses for which continuity is desired. This reflects the two other instances, in which continuity was requested for approximately half of the courses based on the lecturers' input.

To quantify the importance of the subobjectives of the two stakeholder groups and derive the objective weights~$\alpha_1,\ldots,\alpha_{19}$ used in our model, two surveys were conducted at the TUM Campus Straubing in February~2025, one among students and one among lecturers; see Appendix~\ref{UTT:app:survey} for details about these surveys.

\begin{table}[htbp]
\centering
\caption{Summary of instances}
\label{UTT:table:instances}
\begin{tabular}{@{}lccc@{}}
\toprule
 & \textbf{Instance~1} & \textbf{Instance~2} & \textbf{Instance~3} \\
\midrule
Semester                                  & WS 2024/2025 & SS 2025       & WS 2025/2026 \\
Number of courses                         & 181          & 162           & 175          \\
Number of curricula                       & 20           & 20            & 21           \\
Number of modules                         & 92           & 89            & 95           \\
Number of lecturers                       & 93           & 80            & 95           \\
Number of rooms                           & 27           & 27            & 27           \\
Number of periods                         & 30           & 30            & 30           \\
Average number of courses per period         & 6.03         & 5.40           & 5.83         \\
\makecell[l]{Number of courses for which continuity\\of teaching periods is desired}      & 90           & 80            & 82           \\
Average number of courses per curriculum     & 40.15        & 35.30          & 33.67        \\
\bottomrule
\end{tabular}
\end{table}

\subsection{Results}\label{UTT:sec:results}

In this section, we analyze the solution sets obtained for the three instances introduced in Section~\ref{UTT:subsec:instances} under different bounds~$\theta$ on the number of allowed perturbations. For each value~$\theta\in\Theta$, we compute up to ten solutions using the algorithm described in Section~\ref{UTT:subsec:algorithm}, namely two anchor-point solutions and eight additional $\varepsilon$-constraint solutions. Although the algorithm is designed to generate efficient solutions, some computed solutions may still be dominated by other solutions obtained for the same instance and the same perturbation bound~$\theta$ because the individual optimization runs are terminated when a MIP gap of 5\% or a time limit of eight hours is reached. The detailed results of all individual optimization runs, including computation times and MIP gaps, are reported in Tables~\ref{UTT:tab:optimization_runs}, \ref{UTT:tab:optimization_runs2}, and~\ref{UTT:tab:optimization_runs3} in the Appendix. The minimum required number~$\omega$ of perturbations for which a feasible timetable exists is 13 for Instance~1, eight for Instance~2, and 16 for Instance~3.

Figures~\ref{UTT:fig:solution_winter2425}, \ref{UTT:fig:solution_sommer_25}, and~\ref{UTT:fig:solution_winter_2526} show, for the three instances, the obtained images that are not dominated by any other obtained image for the same instance and the same perturbation bound~$\theta$ in the objective space of the lecturers' and students' objectives. The axes of all three plots are truncated to the relevant ranges containing these images.

Tables~\ref{UTT:tab:hypervolumes_inst1}, \ref{UTT:tab:hypervolumes_inst2}, and~\ref{UTT:tab:hypervolumes_inst3} summarize, for each instance and each bound on the number of allowed perturbations, the hypervolume of the corresponding set of images that are not dominated by any other obtained image and the objective values of the anchor points. For each instance, the hypervolumes are computed with respect to a common reference point, namely the Nadir point associated with the union of all image sets obtained for that instance under the considered perturbation bounds. Specifically, for the students' anchor point, we report the students' objective value, while for the lecturers' anchor point, we report the lecturers' objective value. These values represent the best objective values obtained for the respective stakeholder groups under the corresponding perturbation bound. 

We remark that the reported hypervolume values should only be compared across image sets obtained for the same instance. They are not directly comparable across instances, since different instance data can lead to substantially different objective-value ranges. For example, in Instance~1, the lecturers' anchor points for 13, 15, and 17 allowed perturbations have relatively large students' objective values. This is because, with a low number of allowed perturbations and due to the specific time period assignments from the previous year, one course is forced to be scheduled in a Friday evening period, which is particularly heavily penalized in the students' objective.

\subsubsection{Discussion}\label{UTT:subsec:discussion}
\begin{figure}[tbp]
\centering
\begin{tikzpicture}[trim axis group left,trim axis group right]
  \node[rotate=90, anchor=south, yshift=-8.5cm, font=\small] at (-9.6,-2) {Students' objective value};
  \pgfplotsset{
  nadir/.style={only marks, mark=*, mark size=2.2pt, color=black}, anchorpt/.style={only marks, mark=square*, mark size=2.0pt, thick}}
  \begin{groupplot}[
    group style={group size=1 by 2, vertical sep=18pt},
    width=\linewidth,
    xmin=200, 
    xmax=1100,                         
    tick label style={/pgf/number format/use comma},
    legend cell align=left
  ]

\nextgroupplot[
      ymin=30000, 
      ymax=35000,
      axis x line=top,    
      width=\linewidth,
      height=0.35\linewidth,
      grid=both,
      xticklabels=\empty,        
      xlabel={},
      tick style={draw=none},
      axis line style={-},   
      legend style={
      at={(0.95,0.85)}, 
      anchor=north east, 
      fill=white, 
      draw=black, 
      fill opacity=1,
      /tikz/overlay
      },
      clip=true,
      axis on top,
      scaled ticks=false,                        
      tick label style={/pgf/number format/1000 sep={}}  
    ]

\addplot[nadir] coordinates {(1013.333333,34100.91333)};
    \node[font=\small] at (1013.333333,34500) {Nadir point};
    
\addplot[only marks, mark=*, color=blue] coordinates {
    (699,10017.46333) 
    (605,10206.71333)
    (569,10325.11333) 
    (549.6666667,10447.97333) 
    (536.6666667,10770.37333)
    (491,11902.85323) 
    (485.6666667,12010.95333)
    (483,34100.91333) 
  };
  \addplot[anchorpt, color=blue, forget plot] coordinates {
  (483,34100.91333)
  (699,10017.46333)};

  \addplot[only marks, mark=*, color=red] coordinates {
    (875.6666667,9471.883333)
    (744.3333333,9577.693333) 
    (638,9724.833333)
    (564,9866.793333) 
    (498.6666647,10165.42333) 
    (415.6666667,11903.63295)
    (411.6666667,34042.59314) 
  };
  \addplot[anchorpt, color=red, forget plot] coordinates {
  (411.6666667,34042.59314)
  (875.6666667,9471.883333)};

  \addplot[only marks, mark=*, color=green!70!black] coordinates {
    (808,9242.563333) 
    (686.6665967,9349.433268)
    (579,9458.873333) 
    (521.3333333,9547.253333) 
    (466.3333333,9929.793333)
    (382.6666667,11667.93333) 
    (378.6666667,34001.04693) 
  };
  \addplot[anchorpt, color=green!70!black, forget plot] coordinates {
  (378.6666667,34001.04693)
  (808,9242.563333)};

  \addplot[only marks, mark=*, color=orange] coordinates {
    (749,9026.363333)
    (610.6666237,9126.053332)
    (551.3332987,9222.213264)
    (494.3333333,9395.223333)
    (447.3333333,9571.833282)
    (444.6666667,9571.833333)
    (402.3333324,9866.113217)
    (382.3333333,10242.79333)
    (370.6666667,10577.95319)
    (349.3333333,11344.85333)
  };
  \addplot[anchorpt, color=orange, forget plot] coordinates {
  (349.3333333,11344.85333)
  (749,9026.363333)};

  \addplot[only marks, mark=*, color=purple] coordinates {
    (874.3333333,8747.583333)
    (726.3333167,8811.523312)
    (624.6666486,8875.393333)
    (550.6666082,8971.093333)
    (458.6666667,9243.183333)
    (455.3333333,9269.963333)
    (389,9592.673333)
    (367.6666667,9875.013333)
    (343.3333333,10262.79321)
    (323.9999989,10961.49308)
  };
  \addplot[anchorpt, color=purple, forget plot] coordinates {
  (323.9999989,10961.49308)
  (874.3333333,8747.583333)};

  \addplot[only marks, mark=*, color=brown] coordinates {
    (1013.333333,8503.923333)
    (740.6665961,8546.463333)
    (661.9999786,8595.023333)
    (553.9999979,8682.583333)
    (463,9046.953333)
    (376.3333333,9446.773333)
    (345,9805.313333)
    (326,10262.84333)
    (310.6666667,10997.72333)
  };
  \addplot[anchorpt, color=brown, forget plot] coordinates {
  (310.6666667,10997.72333)
  (1013.333333,8503.923333)};

\nextgroupplot[
    ymin=8000, 
    ymax=12500,
    xmin=200,
    axis x line=bottom,
    axis line style={-},
    width=\linewidth,
    height=0.65\linewidth,
    xlabel={Lecturers' objective value},
    tick style={draw=none},
    grid=both,
    legend style={
      at={(0.95,1.4)}, 
      anchor=north east, 
      fill=white, 
      draw=black, 
      fill opacity=1
      },
    scaled ticks=false,                        
    tick label style={/pgf/number format/1000 sep={}},  ]  
  \pgfplotsset{
  nadir/.style={only marks, mark=*, mark size=2.2pt, color=black}, anchorpt/.style={only marks, mark=square*, mark size=2.0pt, thick}}

  \addplot[only marks, mark=*, color=blue] coordinates {
    (699,10017.46333) 
    (605,10206.71333)
    (569,10325.11333) 
    (549.6666667,10447.97333) 
    (536.6666667,10770.37333)
    (491,11902.85323) 
    (485.6666667,12010.95333)
    (483,34100.91333) 
  }; \addlegendentry{13 perturbations}
  \addplot[anchorpt, color=blue, forget plot] coordinates {
  (483,34100.91333)
  (699,10017.46333)};

  \addplot[only marks, mark=*, color=red] coordinates {
    (875.6666667,9471.883333)
    (744.3333333,9577.693333) 
    (638,9724.833333)
    (564,9866.793333) 
    (498.6666647,10165.42333) 
    (415.6666667,11903.63295)
    (411.6666667,34042.59314) 
  }; \addlegendentry{15 perturbations}
  \addplot[anchorpt, color=red, forget plot] coordinates {
  (411.6666667,34042.59314) (875.6666667,9471.883333)};

  \addplot[only marks, mark=*, color=green!70!black] coordinates {
    (808,9242.563333) 
    (686.6665967,9349.433268)
    (579,9458.873333) 
    (521.3333333,9547.253333) 
    (466.3333333,9929.793333)
    (382.6666667,11667.93333) 
    (378.6666667,34001.04693) 
  }; \addlegendentry{17 perturbations}
  \addplot[anchorpt, color=green!70!black, forget plot] coordinates {
  (378.6666667,34001.04693) (808,9242.563333)};

  \addplot[only marks, mark=*, color=orange] coordinates {
    (749,9026.363333)
    (610.6666237,9126.053332)
    (551.3332987,9222.213264)
    (494.3333333,9395.223333)
    (447.3333333,9571.833282)
    (444.6666667,9571.833333)
    (402.3333324,9866.113217)
    (382.3333333,10242.79333)
    (370.6666667,10577.95319)
    (349.3333333,11344.85333)
  }; \addlegendentry{20 perturbations}
  \addplot[anchorpt, color=orange, forget plot] coordinates {
  (349.3333333,11344.85333) (749,9026.363333)};

  \addplot[only marks, mark=*, color=purple] coordinates {
    (874.3333333,8747.583333)
    (726.3333167,8811.523312)
    (624.6666486,8875.393333)
    (550.6666082,8971.093333)
    (458.6666667,9243.183333)
    (455.3333333,9269.963333)
    (389,9592.673333)
    (367.6666667,9875.013333)
    (343.3333333,10262.79321)
    (323.9999989,10961.49308)
  }; \addlegendentry{25 perturbations}
  \addplot[anchorpt, color=purple, forget plot] coordinates {
  (323.9999989,10961.49308) (874.3333333,8747.583333)};

  \addplot[only marks, mark=*, color=brown] coordinates {
    (1013.333333,8503.923333)
    (740.6665961,8546.463333)
    (661.9999786,8595.023333)
    (553.9999979,8682.583333)
    (463,9046.953333)
    (376.3333333,9446.773333)
    (345,9805.313333)
    (326,10262.84333)
    (310.6666667,10997.72333)
  }; \addlegendentry{30 perturbations}
  \addplot[anchorpt, color=brown, forget plot] coordinates {
  (310.6666667,10997.72333) (1013.333333,8503.923333)};

    \addlegendimage{only marks, mark=square, mark size=2.0pt, draw=black, fill=none}
    \addlegendentry{Anchor points}        

  \end{groupplot}

\draw[decorate,decoration={zigzag,amplitude=1pt,segment length=3pt}]
    ($(group c1r1.south west)+(0.02cm,0)$) -- ($(group c1r2.north west)+(0.02cm,0)$);
  \draw[decorate,decoration={zigzag,amplitude=1pt,segment length=3pt}]
    ($(group c1r1.south east)+(-0.02cm,0)$) -- ($(group c1r2.north east)+(-0.02cm,0)$);

\end{tikzpicture}
\caption{Computed images for Instance~1, grouped by perturbation bound.}
\label{UTT:fig:solution_winter2425}
\end{figure}
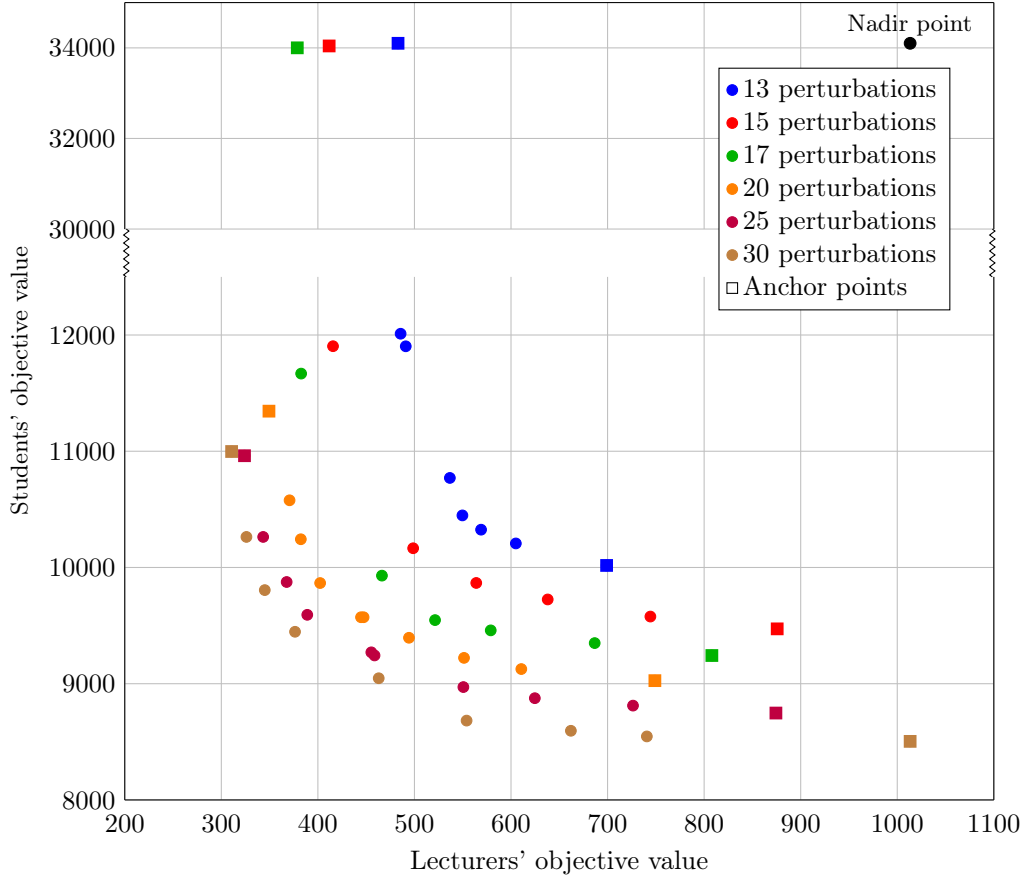

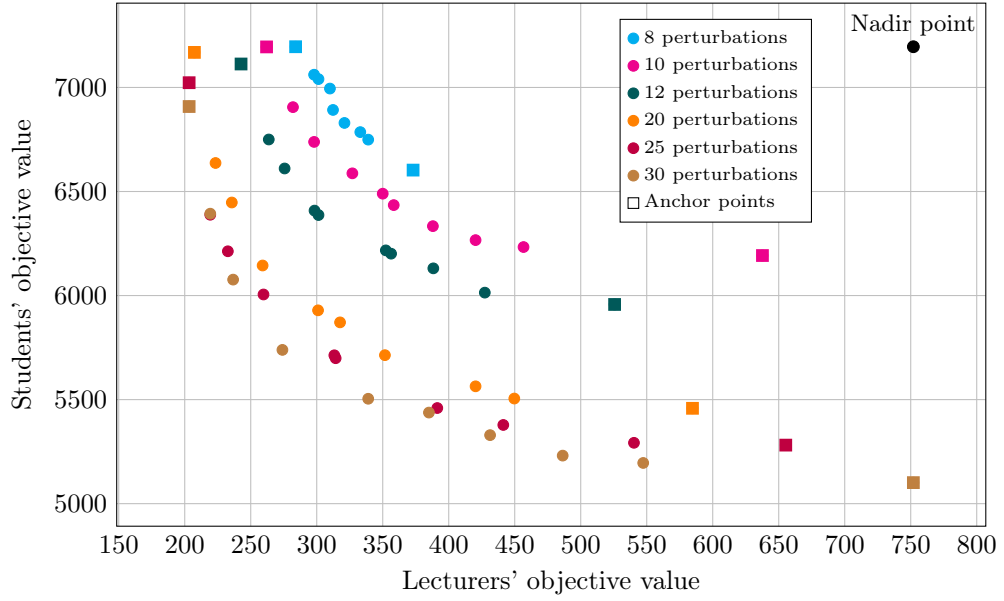
\begin{figure}[tbp]
    \begin{tikzpicture}
      \pgfplotsset{
      nadir/.style={only marks, mark=*, mark size=2.2pt, color=black}, anchorpt/.style={only marks, mark=square*, mark size=2.0pt, thick}}
        \begin{axis}[
            grid=both,
            xlabel={Lecturers' objective value},
            ylabel={Students' objective value},
            width=\linewidth,
            height=0.65\linewidth,
            legend style={at={(0.8,0.97)},anchor=north east, font=\footnotesize},
            legend cell align={left},
            scaled ticks=false,
            tick label style={/pgf/number format/1000 sep={}},
            tick style={draw=none},
            ]
        
\addplot[nadir, forget plot] coordinates {(752,7195.69)};
        \node[font=\small] at (752,7300) {Nadir point};
    
\addplot[
          only marks,
          mark=*,
          color=cyan
        ] coordinates {
          (284,7195.69)
          (298,7061.29)
          (301.3333333,7040.14)
          (309.9999907,6994.679978)
          (312.3333333,6891.98)
(321,6829.48)
          (333,6785.41)
          (339,6749.28)
          (373,6602.9)
        };
        \addlegendentry{8 perturbations}
        \addplot[anchorpt, color=cyan, forget plot] coordinates {
        (284,7195.69)
        (373,6602.9)};

\addplot[
          only marks,
          mark=*,
          color=magenta
        ] coordinates {
          (637.6666667,6192.37)
          (456.6666667,6233.14)
          (420.3333333,6266.43)
          (387.9999795,6333.549964)
          (358.3333333,6434.68)
          (349.9999953,6489.55003)
          (327,6587.299997)
          (298,6738.329974)
          (282,6905.399979)
          (262,7194.83)
        };
        \addlegendentry{10 perturbations}
        \addplot[anchorpt, color=magenta, forget plot] coordinates {
        (262,7194.83)
        (637.6666667,6192.37)};

\addplot[
          only marks,
          mark=*,
          color=teal!70!black
        ] coordinates {
          (525.6666667,5957.19)
          (427.3333333,6014.03)
          (388.3333333,6130.75)
          (356.3333333,6201.22)
          (352.3333333,6217.14)
          (301.3333333,6386.9)
          (298.3333333,6408.22)
          (275.6666667,6610.679995)
          (263.6666667,6749.98)
          (242.6666667,7112.719918)
        };
        \addlegendentry{12 perturbations}
        \addplot[anchorpt, color=teal!70!black, forget plot] coordinates {
        (242.6666667,7112.719918)
        (525.6666667,5957.19)};

\addplot[
          only marks,
          mark=*,
          color=orange
        ] coordinates {
          (584.6666667,5458.19)
          (449.6666667,5505.1)
          (420.3333329,5563.88)
          (351.6666667,5713.56)
          (317.6666667,5871.16)
          (300.9999931,5929.04)
          (259,6144.46)
          (235.6666667,6447.18)
          (223.3333333,6637.01)
          (207.3333333,7168.50993)
        };
        \addlegendentry{20 perturbations}
        \addplot[anchorpt, color=orange, forget plot] coordinates {
        (207.3333333,7168.50993)
        (584.6666667,5458.19)};

\addplot[
          only marks,
          mark=*,
          color=purple
        ] coordinates {
          (655.3333333,5281.31)
          (540.3333229,5292.4)
          (441.3333176,5378.19)
          (391.3333249,5459.52)
          (314.3333333,5698.9)
          (313.3333333,5712.97)
          (259.6666667,6005.16)
          (232.6666667,6212.8)
          (219.3333333,6388.29)
          (203.3333333,7022.85)
        };
        \addlegendentry{25 perturbations}
        \addplot[anchorpt, color=purple, forget plot] coordinates {
        (203.3333333,7022.85)
        (655.3333333,5281.31)};
        
\addplot[
          only marks,
          mark=*,
          color=brown
        ] coordinates {
          (752,5101.12)
          (547.3333011,5195.61)
          (486.3333333,5230.86)
          (431.3333333,5329.279996)
          (385,5437.699952)
          (338.9999985,5504.09)
          (274,5738.9)
          (236.6666667,6076.58)
          (219.3333333,6393.47)
          (203.3333333,6908.379963)
        };
        \addlegendentry{30 perturbations}
        \addplot[anchorpt, color=brown, forget plot] coordinates {
        (203.3333333,6908.379963)
        (752,5101.12)};
        
        \addlegendimage{only marks, mark=square, mark size=2.0pt, draw=black, fill=none}
        \addlegendentry{Anchor points}        
        
      \end{axis}
    \end{tikzpicture}
    \caption{Computed images for Instance~2, grouped by perturbation bound.}
\label{UTT:fig:solution_sommer_25}
\end{figure}

\begin{figure}[tbp]
    \begin{tikzpicture}
      \pgfplotsset{
      nadir/.style={only marks, mark=*, mark size=2.2pt, color=black}, anchorpt/.style={only marks, mark=square*, mark size=2.0pt, thick}}
        \begin{axis}[
            grid=both,
            xlabel={Lecturers' objective value},
            ylabel={Students' objective value},
            width=\linewidth,
            height=0.65\linewidth,
            legend style={                
            at={(0.5,-0.18)},
            anchor=north,
            legend columns=2,
            /tikz/every even column/.append style={column sep=0.5cm},
            },
            legend cell align={left},
            scaled ticks=false,
            tick label style={/pgf/number format/1000 sep={}},
            tick style={draw=none},
            ymin=47700,
            ymax=50450,
]
\addplot[nadir, forget plot] coordinates 
        {(426,50181)};
        \node[font=\small] at (426,50281) {Nadir point};

\addplot[only marks, mark=*, color=violet] coordinates {
        (403, 49578)
        (417, 48939)
        (413, 49037)
        (409, 49196)
        (405, 49333)
        (426, 48839)
        };
        \addlegendentry{16 perturbations}
        \addplot[anchorpt, color=violet, forget plot] coordinates {(403, 49578)(426, 48839)};
        
\addplot[only marks, mark=*, color=olive] coordinates {
        (393, 48979)
        (387, 49217)
        (383, 49393)
        (377, 49537)
        (397, 48803)
        };
        \addlegendentry{18 perturbations}
        \addplot[anchorpt, color=olive, forget plot] coordinates {(377, 49940)(397, 48803)};
        
\addplot[only marks, mark=*, color=orange] coordinates {
        (359, 49533)
        (363, 49437)
        (376, 48927)
        (368, 49170)
        (364, 49373)
        (386, 48646)
        };
        \addlegendentry{20 perturbations}
        \addplot[anchorpt, color=orange, forget plot] coordinates {(359, 50055)(386, 48646)};
        
\addplot[only marks, mark=*, color=purple] coordinates {
        (335, 49213)
        (344, 48748)
        (338, 49034)
        (332, 49288)
        (332, 49288)
        (359, 48406)
        };
        \addlegendentry{25 perturbations}
        \addplot[anchorpt, color=purple, forget plot] coordinates {(332, 50181)(359, 48406)};
        
\addplot[only marks, mark=*, color=brown] coordinates {
        (334, 48834)
        (348, 48453)
        (336, 48697)
        (326, 48952)
        (324, 49169)
        (360, 48209)
        };
        \addlegendentry{30 perturbations}
        \addplot[anchorpt, color=brown, forget plot] coordinates {(324, 49655)(360, 48209)};
        
\addplot[only marks, mark=*, color=black] coordinates {
        (320, 49233)
        (324, 48743)
        (328, 48624)
        (332, 48542)
        (336, 48364)
        (344, 48108)
        (359, 47887)
        };
        \addlegendentry{infinite perturbations}
        \addplot[anchorpt, color=black, forget plot] coordinates {(320, 49233)(359, 47887)};

        \addlegendimage{only marks, mark=square, mark size=2.0pt, draw=black, fill=none}
        \addlegendentry{Anchor points}        
        
      \end{axis}
    \end{tikzpicture}
    \caption{Computed images for Instance~3, grouped by perturbation bound.}\label{UTT:fig:solution_winter_2526}
\end{figure}
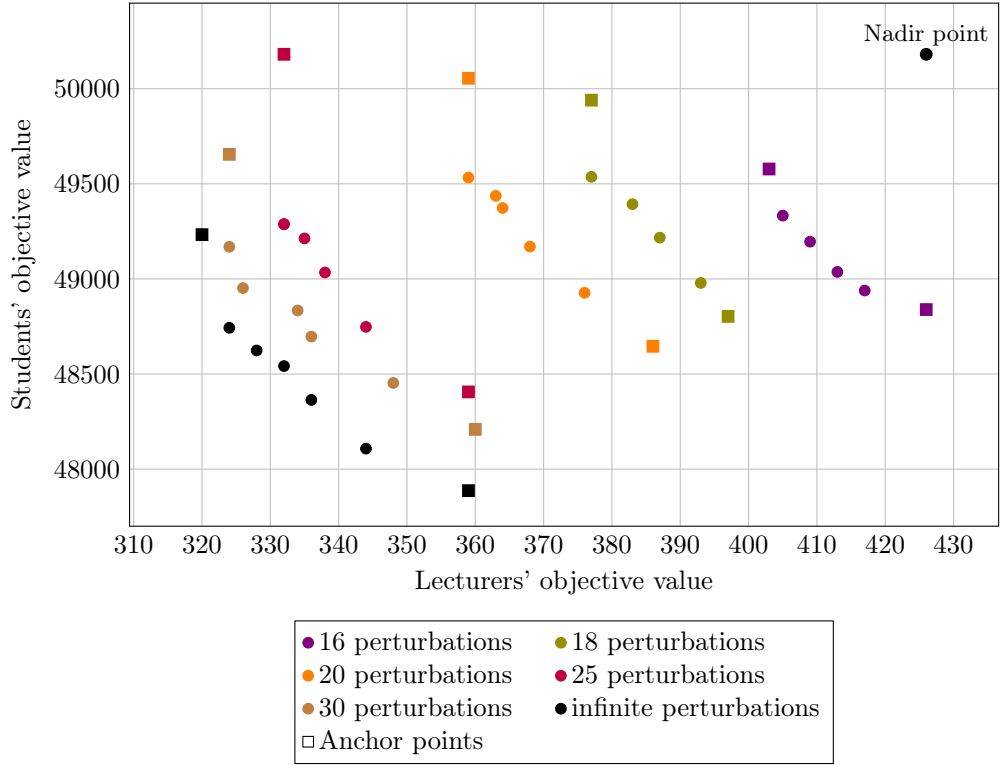

The results reveal a clear trade-off between the lecturers' and students' objectives. As illustrated by the image sets in Figures~\ref{UTT:fig:solution_winter2425}, \ref{UTT:fig:solution_sommer_25}, and~\ref{UTT:fig:solution_winter_2526}, improvements in one objective typically come at the expense of the other, indicating that the two objectives are inherently conflicting. This behavior is consistent across all instances and perturbation bounds considered.

\begin{figure}[htbp]
\centering
\begin{tikzpicture}
\begin{axis}[
    width=8cm, height=8cm,
    enlargelimits=false,
    axis on top,
    colormap/viridis,
    colorbar,
    colorbar style={ylabel={Number of requests}},
    point meta min=0, point meta max=191,
    xlabel={Day}, ylabel={Time period},
    xtick={0,1,2,3,4},
    xticklabels={Mon,Tue,Wed,Thu,Fri},
    ytick={0,1,2,3,4,5},
    yticklabels={8--10,10--12,12--14,14--16,16--18,18--20},
    y dir=reverse,
    tick style={draw=none},
    xticklabel style={anchor=north},
    yticklabel style={anchor=east},
]
\addplot[
    matrix plot*,
    mesh/cols=5, mesh/rows=6,
    point meta=explicit,
    visualization depends on={\thisrow{C} \as \val},
    nodes near coords,
    every node near coord/.append style={
        font=\small,
        /utils/exec={\pgfmathparse{int(\val>95)}\ifnum\pgfmathresult=1\def\labelcolor{black}\else\def\labelcolor{white}\fi},
text/.expanded=\labelcolor,
    },
    nodes near coords={\pgfmathprintnumber[precision=0]{\val}},
] table [meta=C] {
x y C
0 0 71
1 0 74
2 0 106
3 0 97
4 0 49
0 1 112
1 1 140
2 1 165
3 1 191
4 1 89
0 2 97
1 2 136
2 2 113
3 2 169
4 2 67
0 3 116
1 3 135
2 3 95
3 3 131
4 3 17
0 4 48
1 4 64
2 4 60
3 4 61
4 4 4
0 5 0
1 5 17
2 5 24
3 5 7
4 5 1
};
\end{axis}
\end{tikzpicture}
\caption{Number of times each time period was requested as a preferred period, aggregated over the three instances.}
\label{fig:preferred_periods_heatmap}
\end{figure}
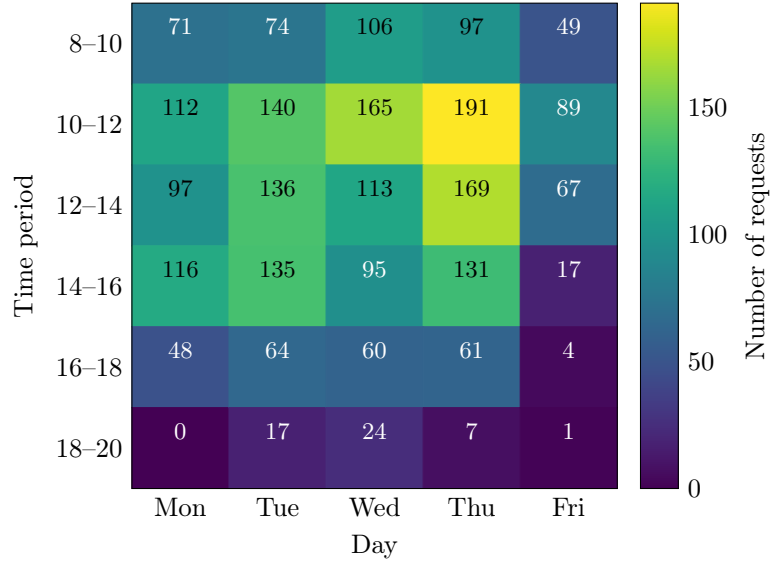

This conflict can be explained by analyzing the structure of the objectives and the instance data. The periods preferred by lecturers are strongly concentrated on Monday to Thursday between 10 a.m. and 4 p.m. (see Figure~\ref{fig:preferred_periods_heatmap}), so reducing the value of lecturers' subobjective~$z_1$ concentrates many courses of a curriculum within the same small set of periods. Since curricula typically contain far more courses than there are such preferred periods (see Table~\ref{UTT:table:instances}), this creates strong pressure toward overlaps within curricula and hence toward a larger value of the students' subobjective~$z_{16}$. A second structural source of conflict concerns consecutive scheduling. Lecturers' subobjective~$z_3$ penalizes gaps between the teaching periods of lecturers who prefer to teach in consecutive periods, whereas the students' subobjectives~$z_{13}$,~$z_{14}$, and~$z_{15}$ penalize consecutive scheduling of lectures or exercise classes within a curriculum. In particular,~$z_{13}$ penalizes two consecutive mandatory lectures or exercise classes of the same curriculum. For such a pair taught by a lecturer who prefers consecutive teaching,~$z_3$ and~$z_{13}$ are naturally conflicting. Scheduling the two courses consecutively improves~$z_3$ but worsens~$z_{13}$. Such pairs are common in our instances: Among all unordered pairs of mandatory lectures or exercise classes within the same curriculum that share at least one lecturer (206 in total across all instances), a shared lecturer prefers consecutive scheduling in 97 cases (47\%).

Moreover, smaller numbers of allowed perturbations limit the attainable timetable quality for both stakeholder groups, highlighting a trade-off between continuity across semesters and timetable quality. Across all instances, relaxing the bound on the number of allowed perturbations substantially improves the achievable objective values for both lecturers and students. This is reflected in the monotonic increase of the hypervolumes. For Instance~1, the hypervolume increases from 125.64~($\times 10^5$) for 13 allowed perturbations to 177.05~($\times 10^5$) for 30 allowed perturbations; see Table~\ref{UTT:tab:hypervolumes_inst1}. Similar trends can be observed for Instance~2 (Table~\ref{UTT:tab:hypervolumes_inst2}) and Instance~3 (Table~\ref{UTT:tab:hypervolumes_inst3}). However, restricting the number of allowed perturbations also substantially reduces the time required per optimization run, as shown in Tables~\ref{UTT:tab:optimization_runs}--\ref{UTT:tab:optimization_runs3}.

Comparing the anchor-point objective values for different perturbation bounds further supports this observation. For both stakeholder groups, relaxing the perturbation bound initially yields substantial improvements in the corresponding anchor-point objective values. The students' objective value at the students' anchor point decreases monotonically as the perturbation bound increases. In contrast, the lecturers' objective value at the lecturers' anchor point ceases to improve substantially beyond certain perturbation bounds in some instances, indicating a saturation effect. For Instance~2, the lecturers' objective value at the lecturers' anchor point remains constant from 25 allowed perturbations onward, and for Instance~3, it improves only marginally from 25 allowed perturbations to the unrestricted case. This suggests that, beyond a certain point, further deviations from the time period assignments in the corresponding semester of the previous academic year do not yield additional benefits for lecturers. The students' objective, on the other hand, continues to improve as more perturbations are allowed.

\begin{table}[tbp]
\caption{Hypervolumes and anchor-point objective values for Instance~1}
\label{UTT:tab:hypervolumes_inst1}
\centering
\renewcommand{\arraystretch}{1.15}
\begin{tabular*}{\textwidth}{@{\extracolsep{\fill}}lrrrrrr@{}}
\toprule
\textbf{Perturbation bound}
    & 13 & 15 & 17 & 20 & 25 & 30 \\
\midrule
Hypervolume ($\times 10^5$)
    & 125.64 & 144.03 & 153.99 & 164.49 & 172.31 & 177.05 \\
\makecell[l]{Lecturers' objective\\at lecturers' anchor point}
    & 483 & 412 & 379 & 349 & 324 & 311 \\
\makecell[l]{Students' objective\\at students' anchor point}
    & 10{,}017 & 9{,}472 & 9{,}243 & 9{,}026 & 8{,}748 & 8{,}504 \\
\bottomrule
\end{tabular*}
\end{table}

\begin{table}[tbp]
\caption{Hypervolumes and anchor-point objective values for Instance~2}
\label{UTT:tab:hypervolumes_inst2}
\centering
\renewcommand{\arraystretch}{1.15}
\begin{tabular*}{\textwidth}{@{\extracolsep{\fill}}lrrrrrr@{}}
\toprule
\textbf{Perturbation bound}
    & 8 & 10 & 12 & 20 & 25 & 30 \\
\midrule
Hypervolume ($\times 10^5$)
    & 2.52 & 4.11 & 5.38 & 8.04 & 8.91 & 9.46 \\
\makecell[l]{Lecturers' objective\\at lecturers' anchor point}
    & 284 & 262 & 242 & 207 & 203 & 203 \\
\makecell[l]{Students' objective\\at students' anchor point}
    & 6{,}602 & 6{,}192 & 5{,}957 & 5{,}458 & 5{,}281 & 5{,}101 \\
\bottomrule
\end{tabular*}
\end{table}

\begin{table}[tbp]
\caption{Hypervolumes and anchor-point objective values for Instance~3}
\label{UTT:tab:hypervolumes_inst3}
\centering
\renewcommand{\arraystretch}{1.15}
\begin{tabular*}{\textwidth}{@{\extracolsep{\fill}}lrrrrrr@{}}
\toprule
\textbf{Perturbation bound}
    & 16 & 18 & 20 & 25 & 30 & $+\infty$ \\
\midrule
Hypervolume ($\times 10^5$)
    & 0.24 & 0.58 & 0.89 & 1.53 & 1.83 & 2.21 \\
\makecell[l]{Lecturers' objective\\at lecturers' anchor point}
    & 403 & 377 & 359 & 332 & 324 & 320 \\
\makecell[l]{Students' objective\\at students' anchor point}
    & 48{,}839 & 48{,}803 & 48{,}646 & 48{,}406 & 48{,}209 & 47{,}887 \\
\bottomrule
\end{tabular*}
\end{table}
 \section{Conclusion}\label{UTT:sec:conclusion}

In this paper, we studied a curriculum-based university course timetabling problem in which the preferences of students and lecturers must be balanced while maintaining continuity across semesters. We proposed a multi-objective mixed-integer programming model that separates the lecturers' and students' objectives and incorporates continuity by limiting the number of allowed perturbations. Based on this model, we developed a multi-objective solution approach that computes sets of candidate timetables and makes the resulting trade-offs explicit.

The computational results on real-world instances from the TUM Campus Straubing show that the lecturers' and students' objectives are inherently conflicting. Thus, relying on a single weighted-sum objective may hide important trade-offs and implicitly encode value judgments that are not made transparent to decision-makers and stakeholders. By contrast, the proposed multi-objective approach provides planners with a set of alternative timetables and allows them to assess how different prioritizations affect the two stakeholder groups. In this sense, the approach supports a more transparent and informed timetabling process.

The results also show that the number of allowed perturbations is a key decision parameter. Strict continuity requirements can substantially restrict the attainable timetable quality for both stakeholder groups, whereas moderate relaxations can lead to considerable improvements. This suggests that continuity across semesters should not be treated as a fixed requirement in isolation, but rather as a planning dimension that needs to be balanced against timetable quality. Making this trade-off explicit can help universities preserve stable organizational routines while still allowing improvements where they provide substantial benefits.

From a practical perspective, the model was developed in close collaboration with the academic planners of the TUM Campus Straubing, and a single-objective version of the mixed-integer programming model presented in this paper was used to generate timetables that were implemented in practice. In addition, we recently developed a web-based graphical user interface that supports the practical timetabling workflow at the campus. The interface allows planners to enter the required information on courses, curricula, and available rooms, enables lecturers to specify their preferences, stores the resulting data in a database, and allows planners to start the optimization process directly from the interface. This demonstrates that the model captures relevant real-world requirements and can be integrated into operational timetable planning. More broadly, the paper illustrates how optimization can be used not only to compute high-quality timetables, but also to structure decision-making in settings where several stakeholder groups and institutional requirements interact.

Future research could build on this work in several directions. First, the approach could be evaluated on data from other universities with different curriculum structures to further assess its transferability. Second, faster algorithms for generating trade-off solutions could be developed to improve scalability for larger instances. Third, the existing graphical user interface could be extended into an interactive decision-support tool that allows planners to explore computed image sets, compare candidate timetables, and communicate the implications of different trade-offs to stakeholders.
 \section*{Statements and Declarations}

\subsection*{Author Contributions}
All authors contributed to the study conception and design. The mathematical model and the multi-objective solution approach were developed jointly by all authors and implemented by Florian Meier. Florian Meier carried out the data collection and wrote the first draft of the manuscript. All authors contributed to the interpretation of the results, revised subsequent versions of the manuscript, and read and approved the final manuscript.

\subsection*{Funding}
\noindent
No funding was received for conducting this study.

\subsection*{Competing Interests}
\noindent
The authors have no relevant financial or non-financial interests to disclose.

\subsection*{Data availability}
\noindent
\sloppy
The datasets generated and analyzed during the current study are not publicly available due to privacy restrictions, as they contain personal data on lecturers' availabilities and preferences as well as internal course and room information from the TUM Campus Straubing.


\begin{thebibliography}{28}
\providecommand{\natexlab}[1]{#1}
\providecommand{\url}[1]{{#1}}
\providecommand{\urlprefix}{URL }
\providecommand{\doi}[1]{\url{https://doi.org/#1}}
\providecommand{\eprint}[2][]{\url{#2}}
 \bibcommenthead

\bibitem[{Akkan et~al.(2022)Akkan, Gülcü, and Kuş}]{akkan_bi-criteria_2022}
Akkan C, Gülcü A, Kuş Z (2022) Bi-criteria simulated annealing for the
  curriculum-based course timetabling problem with robustness approximation.
  Journal of Scheduling 25(4):477--501. \doi{10.1007/s10951-022-00722-0}

\bibitem[{Babaei et~al.(2018)Babaei, Karimpour, and
  Hadidi}]{babaei_applying_2018}
Babaei H, Karimpour J, Hadidi A (2018) Applying hybrid fuzzy multi-criteria
  decision-making approach to find the best ranking for the soft constraint
  weights of lecturers in {UCTP}. International Journal of Fuzzy Systems
  20(1):62--77. \doi{10.1007/s40815-017-0296-z}

\bibitem[{Bellio et~al.(2012)Bellio, Di Gaspero, and
  Schaerf}]{bellio_design_2012}
Bellio R, Di Gaspero L, Schaerf A (2012) Design and statistical analysis of a
  hybrid local search algorithm for course timetabling. Journal of Scheduling
  15(1):49--61. \doi{10.1007/s10951-011-0224-2}

\bibitem[{Bettinelli et~al.(2015)Bettinelli, Cacchiani, Roberti, and
  Toth}]{bettinelli_overview_2015}
Bettinelli A, Cacchiani V, Roberti R, et~al (2015) An overview of
  curriculum-based course timetabling. TOP 23(2):313--349.
  \doi{10.1007/s11750-015-0366-z}

\bibitem[{Ceschia et~al.(2023)Ceschia, Di~Gaspero, and
  Schaerf}]{ceschia_educational_2023}
Ceschia S, Di~Gaspero L, Schaerf A (2023) Educational timetabling: Problems,
  benchmarks, and state-of-the-art results. European Journal of Operational
  Research 308(1):1--18. \doi{10.1016/j.ejor.2022.07.011}

\bibitem[{Chen et~al.(2021)Chen, Sze, Goh, Sabar, and
  Kendall}]{chen_survey_2021}
Chen MC, Sze SN, Goh SL, et~al (2021) A survey of university course timetabling
  problem: Perspectives, trends and opportunities. {IEEE} Access
  9:106515--106529. \doi{10.1109/ACCESS.2021.3100613}

\bibitem[{Chiarandini et~al.(2006)Chiarandini, Birattari, Socha, and
  Rossi-Doria}]{chiarandini_effective_2006}
Chiarandini M, Birattari M, Socha K, et~al (2006) An effective hybrid algorithm
  for university course timetabling. Journal of Scheduling 9(5):403--432.
  \doi{10.1007/s10951-006-8495-8}

\bibitem[{Datta et~al.(2007)Datta, Deb, and Fonseca}]{Datta2007}
Datta D, Deb K, Fonseca CM (2007) Multi-Objective Evolutionary Algorithm for
  University Class Timetabling Problem, Studies in Computational Intelligence
  ({SCI}), vol~49, Springer, pp 197--236. \doi{10.1007/978-3-540-48584-1_8}

\bibitem[{Davison et~al.(2025)Davison, Kheiri, and
  Zografos}]{davison_modelling_2025}
Davison M, Kheiri A, Zografos KG (2025) Modelling and solving the university
  course timetabling problem with hybrid teaching considerations. Journal of
  Scheduling 28(2):195--215. \doi{10.1007/s10951-024-00817-w}

\bibitem[{Dunke and Nickel(2023)}]{dunke_matheuristic_2023}
Dunke F, Nickel S (2023) A matheuristic for customized multi-level
  multi-criteria university timetabling. Annals of Operations Research
  328(2):1313--1348. \doi{10.1007/s10479-023-05325-2}

\bibitem[{Ehrgott(2005)}]{ehrgott_multicriteria_2005}
Ehrgott M (2005) Multicriteria Optimization. Springer-Verlag,
  \doi{10.1007/3-540-27659-9}

\bibitem[{Gülcü and Akkan(2020)}]{gulcu_robust_2020}
Gülcü A, Akkan C (2020) Robust university course timetabling problem subject
  to single and multiple disruptions. European Journal of Operational Research
  283(2):630--646. \doi{10.1016/j.ejor.2019.11.024}

\bibitem[{Hamacher et~al.(2007)Hamacher, Pedersen, and
  Ruzika}]{hamacher_finding_2007}
Hamacher HW, Pedersen CR, Ruzika S (2007) Finding representative systems for
  discrete bicriterion optimization problems. Operations Research Letters
  35(3):336--344. \doi{10.1016/j.orl.2006.03.019}

\bibitem[{Jat and Yang(2011)}]{jat_hybrid_2011}
Jat SN, Yang S (2011) A hybrid genetic algorithm and tabu search approach for
  post enrolment course timetabling. Journal of Scheduling 14(6):617--637.
  \doi{10.1007/s10951-010-0202-0}

\bibitem[{Kuhn et~al.(2016)Kuhn, Fonseca, Paquete, Ruzika, Duarte, and
  Figueira}]{kuhn2016hypervolume}
Kuhn T, Fonseca CM, Paquete L, et~al (2016) Hypervolume subset selection in two
  dimensions: Formulations and algorithms. Evolutionary Computation
  24(3):411--425. \doi{10.1162/EVCO_a_00157}

\bibitem[{Lach and Lübbecke(2012)}]{lach_curriculum_2012}
Lach G, Lübbecke ME (2012) Curriculum based course timetabling: new solutions
  to {Udine} benchmark instances. Annals of Operations Research
  194(1):255--272. \doi{10.1007/s10479-010-0700-7}

\bibitem[{Lemos et~al.(2020)Lemos, Monteiro, and Lynce}]{hebrard_minimal_2020}
Lemos A, Monteiro PT, Lynce I (2020) Minimal perturbation in university
  timetabling with maximum satisfiability. In: Proceedings of the 13th
  International Conference on Integration of Constraint Programming, Artificial
  Intelligence, and Operations Research ({CPAIOR}), pp 317--333,
  \doi{10.1007/978-3-030-58942-4_21}

\bibitem[{Lemos et~al.(2021)Lemos, Monteiro, and
  Lynce}]{lemos_disruptions_2021}
Lemos A, Monteiro PT, Lynce I (2021) Disruptions in timetables: a case study at
  {Universidade de Lisboa}. Journal of Scheduling 24(1):35--48.
  \doi{10.1007/s10951-020-00666-3}

\bibitem[{Li and Yao(2019)}]{li2019quality}
Li M, Yao X (2019) Quality evaluation of solution sets in multiobjective
  optimisation: {A} survey. {ACM} Computing Surveys 52(2):26:1--26:38.
  \doi{10.1145/3300148}

\bibitem[{Lindahl et~al.(2019)Lindahl, Stidsen, and
  Sørensen}]{lindahl_quality_2019}
Lindahl M, Stidsen T, Sørensen M (2019) Quality recovering of university
  timetables. European Journal of Operational Research 276(2):422--435.
  \doi{10.1016/j.ejor.2019.01.026}

\bibitem[{McCollum et~al.(2010)McCollum, Schaerf, Paechter, McMullan, Lewis,
  Parkes, Di~Gaspero, Qu, and Burke}]{mccollum_setting_2010}
McCollum B, Schaerf A, Paechter B, et~al (2010) Setting the research agenda in
  automated timetabling: {The Second International Timetabling Competition}.
  {INFORMS} Journal on Computing 22(1):120--130. \doi{10.1287/ijoc.1090.0320}

\bibitem[{Mühlenthaler and Wanka(2016)}]{muhlenthaler_fairness_2016}
Mühlenthaler M, Wanka R (2016) Fairness in academic course timetabling. Annals
  of Operations Research 239(1):171--188. \doi{10.1007/s10479-014-1553-2}

\bibitem[{Müller et~al.(2005)Müller, Rudová, and
  Barták}]{burke_minimal_2005}
Müller T, Rudová H, Barták R (2005) Minimal perturbation problem in course
  timetabling. In: Proceedings of the 5th International Conference on Practice
  and Theory of Automated Timetabling ({PATAT}), pp 126--146,
  \doi{10.1007/11593577_8}

\bibitem[{Müller et~al.(2025)Müller, Rudová, and
  Müllerová}]{muller_real-world_2025}
Müller T, Rudová H, Müllerová Z (2025) Real-world university course
  timetabling at the {International Timetabling Competition} 2019. Journal of
  Scheduling 28(2):247--267. \doi{10.1007/s10951-023-00801-w}

\bibitem[{Perzina(2007)}]{perzina_solving_2007}
Perzina R (2007) Solving multicriteria university timetabling problem by a
  self-adaptive genetic algorithm with minimal perturbation. In: Proceedings of
  the 5th {IEEE} International Conference on Information Reuse and Integration
  ({IRI}). IEEE, pp 98--103, \doi{10.1109/IRI.2007.4296604}

\bibitem[{Phillips et~al.(2017)Phillips, Walker, Ehrgott, and
  Ryan}]{phillips_integer_2017}
Phillips AE, Walker CG, Ehrgott M, et~al (2017) Integer programming for minimal
  perturbation problems in university course timetabling. Annals of Operations
  Research 252(2):283--304. \doi{10.1007/s10479-015-2094-z}

\bibitem[{Schimmelpfeng and Helber(2007)}]{schimmelpfeng_application_2007}
Schimmelpfeng K, Helber S (2007) Application of a real-world university-course
  timetabling model solved by integer programming. OR Spectrum 29(4):783--803.
  \doi{10.1007/s00291-006-0074-z}

\bibitem[{Vatandoost et~al.(2026)Vatandoost, Golabchi, Ekhlassi, Rahbar, and
  Von~Buelow}]{vatandoost_ranking_2026}
Vatandoost M, Golabchi M, Ekhlassi A, et~al (2026) Ranking and sorting the
  {Pareto} front in optimal shell structure design. Frontiers of Architectural
  Research 15(2):705--728. \doi{10.1016/j.foar.2025.07.006}

\end{thebibliography}

\newpage
\begin{appendices}

\section{Detailed Optimization Results}\label{UTT:subsec:detailed_opt_results}

This appendix reports the detailed results of all individual optimization runs for the three real-world instances considered. Tables~\ref{UTT:tab:optimization_runs}--\ref{UTT:tab:optimization_runs3} contain the objective values, computation times, and MIP gaps of the optimization runs used to compute the solution sets analyzed in Section~\ref{UTT:sec:computational_results}. For each value~$\theta\in\Theta$ for the allowed number of perturbations, the corresponding solutions are grouped together, as indicated in the leftmost column. Each solution is computed lexicographically by first optimizing either the lecturers' or the students' objective and then the other objective under the objective bound obtained from the first run; see Section~\ref{UTT:subsec:algorithm} for details. For each perturbation bound, the first row corresponds to the lecturers' anchor point, while the last row corresponds to the students' anchor point. The tables also report the upper bound~$\parLecturersBound$ or~$\parStudentsBound$ imposed in the first optimization run used to obtain the corresponding solution, if applicable. For the anchor points, no such upper bound is imposed in the first optimization run, and the order of optimization follows from their definition. The tables include all computed solutions, whereas the figures in Section~\ref{UTT:sec:computational_results} omit images dominated by another image obtained for the same instance and perturbation bound.

\small
\begin{longtable}{lllllllll}

\caption{Summary of optimization runs for Instance~1 (see Figure~\ref{UTT:fig:solution_winter2425}). The column ``Max. pert.'' gives the maximum allowed number of perturbations. ``Lect. Obj.'' and ``Stud. Obj.'' report the lecturers' and students' objective values of the computed solution, respectively. ``Time L [s]'' and ``Gap L [\%]'' report the computation time and MIP gap of the optimization run in which the lecturers' objective is optimized, while ``Time S [s]'' and ``Gap S [\%]'' report the corresponding values for the run in which the students' objective is optimized. The columns ``$\parLecturersBound$'' and ``$\parStudentsBound$'' provide the upper bounds imposed on the lecturers' and students' objectives, respectively, in the first optimization run.}
\label{UTT:tab:optimization_runs}\\

\toprule
\rotatebox{90}{\makebox{Max. pert.}} &
\rotatebox{90}{\makebox{Lect. Obj.}} &
\rotatebox{90}{\makebox{Stud. Obj.}} &
\rotatebox{90}{\makebox{Time L [s]}} &
\rotatebox{90}{\makebox{Gap L [\%]}} &
\rotatebox{90}{\makebox{$\parLecturersBound$}} &
\rotatebox{90}{\makebox{Time S [s]}} &
\rotatebox{90}{\makebox{Gap S [\%]}} &
\rotatebox{90}{\makebox{$\parStudentsBound$}} \rule{0pt}{3.5em} \\
\midrule
\endfirsthead

\multicolumn{7}{c}{\tablename\ \thetable\ -- \textit{continued from previous page}} \\
\toprule
\rotatebox{90}{\makebox{Max. pert.}} &
\rotatebox{90}{\makebox{Lect. Obj.}} &
\rotatebox{90}{\makebox{Stud. Obj.}} &
\rotatebox{90}{\makebox{Time L [s]}} &
\rotatebox{90}{\makebox{Gap L [\%]}} &
\rotatebox{90}{\makebox{$\parLecturersBound$}} &
\rotatebox{90}{\makebox{Time S [s]}} &
\rotatebox{90}{\makebox{Gap S [\%]}} &
\rotatebox{90}{\makebox{$\parStudentsBound$}} \rule{0pt}{3.5em} \\
\midrule
\endhead

\midrule
\multicolumn{7}{r}{\textit{continued on next page}} \\
\endfoot

\bottomrule
\endlastfoot

13 & 483 & 34101 & 7 & 4.62 & \NA & 90 & 0.39 & \NA\\
13 & 491 & 11903 & 12 & 4.89 & 537 & 51 & 4.99 & \NA \\
13 & 486 & 12031 & 19 & 3.16 & 573 & 40 & 4.93 & \NA \\
13 & 486 & 12102 & 22 & 4.74 & 609 & 49 & 4.77 & \NA \\
13 & 486 & 12011 & 19 & 3.36 & 645 & 57 & 4.81 & \NA \\
13 & 537 & 10770 & 33 & 4.82 & \NA & 40 & 4.42 & 16038 \\
13 & 550 & 10448 & 52 & 4.50 & \NA & 30 & 4.18 & 20052 \\
13 & 569 & 10325 & 89 & 4.78 & \NA & 23 & 4.99 & 24066 \\
13 & 605 & 10207 & 136 & 4.86 & \NA & 18 & 4.92 & 28080 \\
13 & 699 & 10017 & 213 & 4.72 & \NA & 13 & 4.01 & \NA \\

\addlinespace

15 & 412& 34043 & 11 & 3.89 & \NA & 539 & 0.02 & \NA \\
15 & 417& 11932 & 44 & 4.88 & 527 & 232 & 4.63 & \NA \\
15 & 416& 11904 & 27 & 4.97 & 605 & 298 & 4.95 & \NA \\
15 & 416& 11932 & 35 & 4.97 & 682 & 114 & 4.83 & \NA \\
15 & 416& 11904 & 44 & 4.97 & 760 & 270 & 4.36 & \NA \\
15 & 499 & 10165 & 638 & 4.97 & \NA & 581 & 4.84 & 15615 \\
15 & 564 & 9867 & 4731 & 4.96 & \NA & 287 & 4.92 & 19710 \\
15 & 638 & 9725 & 7017 & 4.86 & \NA & 331 & 4.87 & 23805 \\
15 & 744 & 9578 & 7061 & 4.96 & \NA & 237 & 4.94 & 27900 \\
15 & 876 & 9472 & 11373 & 4.87 & \NA & 153 & 4.94 & \NA \\

\addlinespace

17 & 379 & 34001 & 16 & 4.05 & \NA & 536 & 0.52 & \NA \\
17 & 388 & 11470 & 39 & 4.64 & 486 & 299 & 4.79 & \NA \\
17 & 383 & 11668 & 29 & 4.97 & 558 & 382 & 4.81 & \NA \\
17 & 383 & 11688 & 42 & 4.97 & 629 & 436 & 4.98 & \NA \\
17 & 383 & 11668 & 40 & 4.88 & 701 & 366 & 4.82 & \NA \\
17 & 466 & 9930 & 805 & 4.91 & \NA & 635 & 4.78 & 15432 \\
17 & 521 & 9547 & 6627 & 4.74 & \NA & 382 & 5.00 & 19559 \\
17 & 579 & 9459 & 8895 & 4.93 & \NA & 480 & 4.95 & 23685 \\
17 & 687 & 9349 & 28800 & 7.02 & \NA & 416 & 4.93 & 27811 \\
17 & 808 & 9243 & 28800 & 9.98 & \NA & 194 & 4.98 & \NA \\

\addlinespace

20 & 349 & 11345 & 19 & 4.77 & \NA & 462 & 4.34 & \NA \\
20 & 447 & 9572 & 8747 & 4.62 & 449 & 1082 & 4.79 & \NA \\
20 & 402 & 9866 & 694 & 4.97 & 516 & 687 & 4.96 & \NA \\
20 & 382 & 10243 & 346 & 4.97 & 582 & 621 & 4.94 & \NA \\
20 & 371 & 10578 & 258 & 4.86 & 649 & 731 & 4.97 & \NA \\
20 & 445 & 9572 & 10017 & 4.89 & \NA & 920 & 4.66 & 9606 \\
20 & 494 & 9395 & 12961 & 4.93 & \NA & 1543 & 4.93 & 9992 \\
20 & 551 & 9222 & 28800 & 8.27 & \NA & 1421 & 4.66 & 10379 \\
20 & 611 & 9126 & 28800 & 12.59 & \NA & 1018 & 4.89 & 10765 \\
20 & 749 & 9026 & 28801 & 24.33 & \NA & 643 & 4.95 & \NA \\

\addlinespace

25 & 324 & 10961 & 46 & 4.94 & \NA & 982 & 4.93 & \NA \\
25 & 455 & 9270 & 28800 & 11.27 & 462 & 13845 & 4.94 & \NA \\
25 & 389 & 9593 & 8210 & 4.95 & 553 & 2521 & 4.99 & \NA \\
25 & 368 & 9875 & 5638 & 4.71 & 645 & 2012 & 4.96 & \NA \\
25 & 343 & 10263 & 643 & 4.95 & 737 & 1331 & 4.46 & \NA \\
25 & 459 & 9243 & 28800 & 12.98 & \NA & 14050 & 4.75 & 9301 \\
25 & 551 & 8971 & 28800 & 18.49 & \NA & 14635 & 4.87 & 9670 \\
25 & 625 & 8875 & 28800 & 26.47 & \NA & 15473 & 4.83 & 10039 \\
25 & 726 & 8812 & 28800 & 34.91 & \NA & 14950 & 4.74 & 10408 \\
25 & 874 & 8748 & 28800 & 45.74 & \NA & 8521 & 4.98 & \NA \\

\addlinespace

30 & 311 & 10998 & 99 & 4.91 & \NA & 10544 & 4.89 & \NA \\
30 & 463 & 9047 & 28800 & 18.79 & 486 & 28800 & 5.40 & \NA \\
30 & 376 & 9447 & 26384 & 4.96 & 603 & 16094 & 4.96 & \NA \\
30 & 345 & 9805 & 8113 & 4.44 & 721 & 12576 & 4.97 & \NA \\
30 & 326 & 10263 & 1113 & 4.91 & 838 & 9307 & 4.74 & \NA \\
30 & 483 & 9051 & 28800 & 21.14 & \NA & 28800 & 6.58 & 9127 \\
30 & 554 & 8683 & 28800 & 24.20 & \NA & 28367 & 4.68 & 9543 \\
30 & 662 & 8595 & 28800 & 37.82 & \NA & 23509 & 5.00 & 9959 \\
30 & 741 & 8546 & 28800 & 44.03 & \NA & 19848 & 5.00 & 10374 \\
30 & 1013 & 8504 & 28800 & 58.15 & \NA & 12074 & 4.62 & \NA \\

\end{longtable}

\small
\begin{longtable}{lllllllll}

\caption{Summary of optimization runs for \textbf{Instance~2} (see~\ Figure~\ref{UTT:fig:solution_sommer_25}). The table is structured analogously to Table~\ref{UTT:tab:optimization_runs} for Instance~1.\label{UTT:tab:optimization_runs2}}\\

\toprule
\rotatebox{90}{\makebox{Max. pert.}} &
\rotatebox{90}{\makebox{Lect. Obj.}} &
\rotatebox{90}{\makebox{Stud. Obj.}} &
\rotatebox{90}{\makebox{Time L [s]}} &
\rotatebox{90}{\makebox{Gap L [\%]}} &
\rotatebox{90}{\makebox{$\parLecturersBound$}} &
\rotatebox{90}{\makebox{Time S [s]}} &
\rotatebox{90}{\makebox{Gap S [\%]}} &
\rotatebox{90}{\makebox{$\parStudentsBound$}} \rule{0pt}{3.5em} \\
\midrule
\endfirsthead

\multicolumn{7}{c}{\tablename\ \thetable\ -- \textit{continued from previous page}} \\
\toprule
\rotatebox{90}{\makebox{Max. pert.}} &
\rotatebox{90}{\makebox{Lect. Obj.}} &
\rotatebox{90}{\makebox{Stud. Obj.}} &
\rotatebox{90}{\makebox{Time L [s]}} &
\rotatebox{90}{\makebox{Gap L [\%]}} &
\rotatebox{90}{\makebox{$\parLecturersBound$}} &
\rotatebox{90}{\makebox{Time S [s]}} &
\rotatebox{90}{\makebox{Gap S [\%]}} &
\rotatebox{90}{\makebox{$\parStudentsBound$}} \rule{0pt}{3.5em} \\
\midrule
\endhead

\midrule
\multicolumn{7}{r}{\textit{continued on next page}} \\
\endfoot

\bottomrule
\endlastfoot

8 & 284 & 7196 & 13 & 2.00 & \NA & 69 & 2.71 & \NA \\
8 & 298 & 7061 & 74 & 3.90 & 306 & 43 & 4.88 & \NA \\
8 & 310 & 6995 & 132 & 5.00 & 321 & 45 & 5.00 & \NA \\
8 & 316 & 6920 & 164 & 4.99 & 336 & 43 & 4.93 & \NA \\
8 & 333 & 6785 & 85 & 4.89 & 351 & 33 & 4.04 & \NA \\
8 & 339 & 6749 & 148 & 4.13 & \NA & 16 & 4.92 & 6751 \\
8 & 321 & 6829 & 109 & 3.74 & \NA & 20 & 4.22 & 6850 \\
8 & 312 & 6892 & 98 & 3.63 & \NA & 21 & 4.57 & 6949 \\
8 & 301 & 7040 & 83 & 4.98 & \NA & 24 & 4.99 & 7047 \\
8 & 373 & 6602 & 139 & 4.65 & \NA & 9 & 4.61 & \NA \\

\addlinespace

10 & 262 & 7195 & 24 & 4.96 & \NA & 266 & 4.56 & \NA \\
10 & 350 & 6490 & 1228 & 4.94 & 356 & 388 & 4.99 & \NA \\
10 & 388 & 6334 & 1159 & 4.99 & 419 & 102 & 4.74 & \NA \\
10 & 420 & 6266 & 1240 & 4.87 & 481 & 147 & 4.93 & \NA \\
10 & 457 & 6233 & 1948 & 3.38 & 544 & 82 & 4.93 & \NA \\
10 & 358 & 6435 & 2010 & 4.93 & \NA & 350 & 4.48 & 6443 \\
10 & 327 & 6587 & 1645 & 4.99 & \NA & 296 & 4.72 & 6610 \\
10 & 298 & 6738 & 1368 & 4.92 & \NA & 321 & 4.26 & 6777 \\
10 & 282 & 6905 & 710 & 4.02 & \NA & 650 & 4.67 & 6944 \\
10 & 638 & 6192 & 1514 & 0.20 & \NA & 34 & 4.03 & \NA \\

\addlinespace

12 & 242 & 7113 & 38 & 4.95 & \NA & 611 & 4.36 & \NA \\
12 & 298 & 6408 & 3799 & 4.98 & 313 & 747 & 4.45 & \NA \\
12 & 352 & 6217 & 7106 & 4.95 & 361 & 932 & 4.88 & \NA \\
12 & 388 & 6131 & 5293 & 4.85 & 408 & 476 & 5.00 & \NA \\
12 & 427 & 6014 & 2669 & 4.82 & 455 & 477 & 4.99 & \NA \\
12 & 356 & 6201 & 6072 & 4.96 & \NA & 485 & 4.91 & 6246 \\
12 & 301 & 6387 & 1649 & 4.98 & \NA & 652 & 4.79 & 6439 \\
12 & 276 & 6611 & 1236 & 4.11 & \NA & 1089 & 4.77 & 6631 \\
12 & 264 & 6750 & 830 & 4.68 & \NA & 767 & 4.48 & 6824 \\
12 & 526 & 5957 & 1449 & 0.49 & \NA & 67 & 4.99 & \NA \\

\addlinespace

20 & 207 & 7169 & 47 & 4.50 & \NA & 1198 & 4.48 & \NA \\
20 & 301 & 5929 & 28800 & 6.72 & 302 & 8021 & 4.95 & \NA \\
20 & 352 & 5714 & 28800 & 13.23 & 365 & 12604 & 4.94 & \NA \\
20 & 420 & 5564 & 28800 & 20.24 & 427 & 10784 & 4.84 & \NA \\
20 & 450 & 5505 & 28800 & 19.21 & 490 & 2560 & 4.90 & \NA \\
20 & 318 & 5871 & 28800 & 10.81 & \NA & 11392 & 4.96 & 5886 \\
20 & 259 & 6144 & 13956 & 4.89 & \NA & 11705 & 5.00 & 6171 \\
20 & 236 & 6447 & 3844 & 4.95 & \NA & 3161 & 4.95 & 6456 \\
20 & 223 & 6637 & 1317 & 4.77 & \NA & 2294 & 4.97 & 6741 \\
20 & 585 & 5458 & 28800 & 36.10 & \NA & 711 & 4.69 & \NA \\

\addlinespace

25 & 203 & 7023 & 57 & 4.92 & \NA & 843 & 4.70 & \NA \\
25 & 314 & 5699 & 28800 & 17.70 & 316 & 25770 & 4.98 & \NA \\
25 & 391 & 5460 & 28800 & 25.74 & 392 & 20234 & 4.99 & \NA \\
25 & 441 & 5378 & 28800 & 32.67 & 467 & 18566 & 4.95 & \NA \\
25 & 540 & 5292 & 28800 & 40.65 & 542 & 15753 & 4.97 & \NA \\
25 & 313 & 5713 & 28800 & 15.85 & \NA & 26819 & 4.97 & 5717 \\
25 & 260 & 6005 & 28221 & 4.88 & \NA & 14050 & 4.86 & 6007 \\
25 & 233 & 6213 & 11408 & 4.96 & \NA & 12863 & 4.56 & 6297 \\
25 & 219 & 6388 & 2246 & 4.86 & \NA & 1602 & 4.95 & 6587 \\
25 & 655 & 5281 & 28800 & 52.00 & \NA & 7174 & 4.88 & \NA \\

\addlinespace

30 & 203 & 6908 & 48 & 4.92 & \NA & 1002 & 4.85 & \NA \\
30 & 339 & 5504 & 28800 & 23.06 & 340 & 28800 & 7.36 & \NA \\
30 & 431 & 5329 & 28800 & 36.21 & 432 & 28800 & 6.95 & \NA \\
30 & 486 & 5231 & 28800 & 42.82 & 523 & 28800 & 5.90 & \NA \\
30 & 547 & 5196 & 28801 & 47.60 & 615 & 28800 & 5.13 & \NA \\
30 & 385 & 5438 & 28800 & 33.68 & \NA & 28800 & 7.76 & 5553 \\
30 & 274 & 5739 & 28800 & 11.44 & \NA & 28800 & 5.51 & 5854 \\
30 & 237 & 6077 & 13363 & 4.93 & \NA & 14912 & 4.99 & 6155 \\
30 & 219 & 6393 & 2138 & 4.83 & \NA & 11237 & 4.21 & 6457 \\
30 & 752 & 5101 & 28800 & 60.32 & \NA & 8442 & 4.96 & \NA \\

\end{longtable}

\small
\begin{longtable}{lllllllll}

\caption{Summary of optimization runs for \textbf{Instance~3} (see~\ Figure~\ref{UTT:fig:solution_winter_2526}). The table is structured analogously to Table~\ref{UTT:tab:optimization_runs} for Instance~1.\label{UTT:tab:optimization_runs3}}\\

\toprule
\rotatebox{90}{\makebox{Max. pert.}} &
\rotatebox{90}{\makebox{Lect. Obj.}} &
\rotatebox{90}{\makebox{Stud. Obj.}} &
\rotatebox{90}{\makebox{Time L [s]}} &
\rotatebox{90}{\makebox{Gap L [\%]}} &
\rotatebox{90}{\makebox{$\parLecturersBound$}} &
\rotatebox{90}{\makebox{Time S [s]}} &
\rotatebox{90}{\makebox{Gap S [\%]}} &
\rotatebox{90}{\makebox{$\parStudentsBound$}} \rule{0pt}{3.5em} \\
\midrule
\endfirsthead

\multicolumn{7}{c}{\tablename\ \thetable\ -- \textit{continued from previous page}} \\
\toprule
\rotatebox{90}{\makebox{Max. pert.}} &
\rotatebox{90}{\makebox{Lect. Obj.}} &
\rotatebox{90}{\makebox{Stud. Obj.}} &
\rotatebox{90}{\makebox{Time L [s]}} &
\rotatebox{90}{\makebox{Gap L [\%]}} &
\rotatebox{90}{\makebox{$\parLecturersBound$}} &
\rotatebox{90}{\makebox{Time S [s]}} &
\rotatebox{90}{\makebox{Gap S [\%]}} &
\rotatebox{90}{\makebox{$\parStudentsBound$}} \rule{0pt}{3.5em} \\
\midrule
\endhead

\midrule
\multicolumn{7}{r}{\textit{continued on next page}} \\
\endfoot

\bottomrule
\endlastfoot

16 & 403 & 49578 & 8 & 0.05 & \NA & 56 & 0.02 & \NA \\
16 & 405 & 49417 & 27 & 0.05 & 409 & 9 & 0.05 & \NA \\
16 & 405 & 49417 & 23 & 0.05 & 413 & 6 & 0.05 & \NA \\
16 & 405 & 49417 & 24 & 0.05 & 417 & 6 & 0.05 & \NA \\
16 & 405 & 49417 & 22 & 0.05 & 420 & 6 & 0.05 & \NA \\
16 & 417 & 48939 & 48 & 0.05 & \NA & 6 & 0.05 & 49024 \\
16 & 413 & 49037 & 46 & 0.05 & \NA & 6 & 0.05 & 49147 \\
16 & 409 & 49196 & 41 & 0.05 & \NA & 6 & 0.05 & 49270 \\
16 & 405 & 49333 & 34 & 0.04 & \NA & 6 & 0.05 & 49393 \\
16 & 426 & 48839 & 81 & 0.04 & \NA & 4 & 0.04 & \NA \\

\addlinespace

18 & 377 & 49940 & 18 & 0.05 & \NA & 21 & 0.04 & \NA \\
18 & 377 & 49815 & 57 & 0.05 & 382 & 12 & 0.05 & \NA \\
18 & 377 & 49599 & 46 & 0.05 & 385 & 12 & 0.05 & \NA \\
18 & 377 & 49590 & 41 & 0.05 & 388 & 10 & 0.05 & \NA \\
18 & 377 & 49570 & 56 & 0.05 & 392 & 10 & 0.05 & \NA \\
18 & 393 & 48979 & 292 & 0.05 & \NA & 6 & 0.05 & 49088 \\
18 & 387 & 49217 & 140 & 0.05 & \NA & 10 & 0.05 & 49277 \\
18 & 383 & 49393 & 98 & 0.05 & \NA & 11 & 0.05 & 49466 \\
18 & 377 & 49537 & 115 & 0.05 & \NA & 9 & 0.05 & 49656 \\
18 & 397 & 48803 & 288 & 0.05 & \NA & 4 & 0.05 & \NA \\

\addlinespace

20 & 359 & 50055 & 23 & 0.05 & \NA & 29 & 0.05 & \NA \\
20 & 359 & 49579 & 42 & 0.05 & 365 & 13 & 0.05 & \NA \\
20 & 359 & 49578 & 35 & 0.05 & 370 & 13 & 0.05 & \NA \\
20 & 359 & 49533 & 30 & 0.05 & 375 & 17 & 0.05 & \NA \\
20 & 363 & 49437 & 67 & 0.05 & 379 & 13 & 0.05 & \NA \\
20 & 376 & 48927 & 331 & 0.05 & \NA & 11 & 0.05 & 48998 \\
20 & 368 & 49170 & 160 & 0.05 & \NA & 11 & 0.05 & 49233 \\
20 & 364 & 49373 & 88 & 0.05 & \NA & 11 & 0.05 & 49468 \\
20 & 359 & 49573 & 56 & 0.05 & \NA & 12 & 0.05 & 49703 \\
20 & 386 & 48646 & 401 & 0.05 & \NA & 6 & 0.05 & \NA \\

\addlinespace

25 & 332 & 50181 & 22 & 0.05 & \NA & 30 & 0.05 & \NA \\
25 & 332 & 49866 & 75 & 0.05 & 339 & 29 & 0.05 & \NA \\
25 & 332 & 49616 & 85 & 0.04 & 343 & 22 & 0.05 & \NA \\
25 & 332 & 49294 & 91 & 0.04 & 348 & 17 & 0.05 & \NA \\
25 & 335 & 49213 & 110 & 0.05 & 352 & 15 & 0.05 & \NA \\
25 & 344 & 48748 & 403 & 0.05 & \NA & 10 & 0.05 & 48850 \\
25 & 338 & 49034 & 188 & 0.05 & \NA & 11 & 0.05 & 49146 \\
25 & 332 & 49288 & 118 & 0.05 & \NA & 14 & 0.05 & 49441 \\
25 & 332 & 49288 & 110 & 0.05 & \NA & 15 & 0.05 & 49737 \\
25 & 359 & 48406 & 480 & 0.05 & \NA & 10 & 0.05 & \NA \\

\addlinespace

30 & 324 & 49655 & 20 & 0.05 & \NA & 54 & 0.05 & \NA \\
30 & 324 & 49275 & 51 & 0.04 & 333 & 23 & 0.05 & \NA \\
30 & 327 & 49081 & 122 & 0.04 & 339 & 21 & 0.05 & \NA \\
30 & 327 & 48989 & 119 & 0.05 & 345 & 19 & 0.05 & \NA \\
30 & 334 & 48834 & 246 & 0.05 & 351 & 21 & 0.05 & \NA \\
30 & 348 & 48453 & 978 & 0.05 & \NA & 11 & 0.05 & 48570 \\
30 & 336 & 48697 & 391 & 0.05 & \NA & 17 & 0.05 & 48811 \\
30 & 326 & 48952 & 332 & 0.05 & \NA & 19 & 0.05 & 49052 \\
30 & 324 & 49169 & 164 & 0.05 & \NA & 21 & 0.05 & 49293 \\
30 & 360 & 48209 & 1544 & 0.05 & \NA & 13 & 0.05 & \NA \\

\addlinespace

inf. & 320 & 49233 & 21 & 0.04 & \NA & 154 & 0.04 & \NA \\
inf. & 326 & 48937 & 247 & 0.05 & 330 & 41 & 0.05 & \NA \\
inf. & 328 & 48765 & 308 & 0.05 & 337 & 37 & 0.05 & \NA \\
inf. & 328 & 48709 & 211 & 0.05 & 343 & 41 & 0.05 & \NA \\
inf. & 328 & 48624 & 262 & 0.05 & 350 & 41 & 0.05 & \NA \\
inf. & 344 & 48108 & 1988 & 0.05 & \NA & 18 & 0.05 & 48223 \\
inf. & 336 & 48364 & 847 & 0.05 & \NA & 24 & 0.05 & 48448 \\
inf. & 332 & 48542 & 834 & 0.05 & \NA & 30 & 0.05 & 48672 \\
inf. & 324 & 48743 & 293 & 0.05 & \NA & 30 & 0.05 & 48897 \\
inf. & 359 & 47887 & 9979 & 0.05 & \NA & 38 & 0.05 & \NA \\

\end{longtable}
\normalsize

\section{Surveys}\label{UTT:app:survey}

The surveys mentioned in Section~\ref{UTT:subsec:instances} were used to quantify the relative importance of the subobjectives and to derive the weights~$\alpha_1,\ldots,\alpha_{19}$. Separate surveys were conducted for students and lecturers. In each survey, the relevant subobjectives were first explained, and respondents were then asked to rank them from highest to lowest priority and assign each subobjective an importance score between zero and ten.

To validate the stated preferences, respondents were additionally presented with pairwise comparisons of potentially conflicting subobjectives. These comparisons were accompanied by concrete scheduling examples, and respondents were asked which alternative they preferred. The rankings, importance scores, and pairwise comparisons were then checked for consistency. If, for a given respondent, the relative order of two subobjectives was inconsistent with their assigned scores or with the stated preference in a corresponding pairwise comparison, the scores provided by this respondent for the affected subobjectives were excluded.

Based on the remaining responses, an average importance score between zero and ten was computed for each subobjective and for the corresponding stakeholder group. To obtain suitable weights for the MIP model, these scores were normalized by the typical magnitudes of the subobjective values. For this purpose, we solved the MIP model with all weights~$\alpha_1,\ldots,\alpha_{19}$ set to one. Since only Instances~1 and~2 were available when this normalization was performed, we used these two instances for the calibration. The resulting subobjective values were very similar across the two instances, and we therefore used their mean values as normalization factors. The weight~$\alpha_i$ of each subobjective was then computed by dividing its survey-based importance score by the corresponding normalization factor. Thus, the weights reflect both the perceived importance of the subobjectives and their relative magnitudes in the model.

The student survey was also used to determine the period-specific penalties~$\parPeriodPenalty(p)$ for each period~$p \in \setPeriods$. Students were asked which periods they preferred to avoid and how inconvenient they considered courses in these periods, choosing among ``Extreme Inconvenience'', ``Significant Inconvenience'', ``Moderate Inconvenience'', and ``Minimal Inconvenience''. The penalty~$\parPeriodPenalty(p)$ was computed as a weighted sum of the number of times period~$p$ was rated as inconvenient, using weights 10, 6, 3, and 1 for the four inconvenience levels, respectively.

In total, 91 students and 60 lecturers participated in the surveys. The resulting weights for both stakeholder groups are summarized in Table~\ref{UTT:tab:weights}.

\begin{table}[htbp]
\centering
\caption{Survey-based weights for the subobjectives of both stakeholder groups.}
\label{UTT:tab:weights}
\begin{tabular}{ccccccccccc}
\toprule
Lecturers & $\alpha_1$ & $\alpha_2$ & $\alpha_3$ & $\alpha_4$ & $\alpha_5$ & $\alpha_6$ & $\alpha_7$ & $\alpha_8$ & $\alpha_9$ & $\alpha_{10}$ \\
\midrule
 & 4 & 14 & 12 & 37 & 11 & 11 & 14 & 14 & 3 & 4 \\
\midrule
Students & $\alpha_{11}$ & $\alpha_{12}$ & $\alpha_{13}$ & $\alpha_{14}$ & $\alpha_{15}$ & $\alpha_{16}$ & $\alpha_{17}$ & $\alpha_{18}$ & $\alpha_{19}$\\
\midrule
 & 95 & 100 & 10 & 23 & 24 & 10 & 50 & 6 & 1 \\
\bottomrule
\end{tabular}
\end{table}

\end{appendices} 
\clearpage

\end{document}